\documentclass[10pt, a4paper]{article}
\usepackage{cmap}
\usepackage[main=russian,english]{babel}
\usepackage[utf8]{inputenc}
\usepackage[text={7in,10in},centering]{geometry}
\usepackage{amsfonts,amsmath,amsxtra,amsthm,amssymb,latexsym, tabularray}
\usepackage{graphicx}
\usepackage{verbatim}
\usepackage{cite, color}
\usepackage{subcaption}
\usepackage{float}
\theoremstyle{definition}
\newtheorem*{definition}{Определение}

\DeclareMathOperator{\sign}{sign}
\DeclareMathOperator{\arsh}{arsh}

\newcommand{\VEC}{\operatorname{Vec}\nolimits}

\newcommand{\dif}{\mathrm{d}}

\newcommand{\be}[1]{\begin{equation}\label{#1}}
\newcommand{\ee}{\end{equation}}

\begin{document}
\parindent=0cm

\title{Лоренцева плоскость анти-де Ситтера}
\author{А.З. Али,\ Ю.Л. Сачков}
\date{\today}

\maketitle

\begin{abstract}
В данной работе исследуется двумерная лоренцева задача на плоскости анти-де Ситтера. Используя методы геометрической теории управления и дифференциальной геометрии, удалось построить ортонормированный репер, вычислить экстремальные траектории, описать множество достижимости, построить оптимальный синтез, описать лоренцево расстояние.
\end{abstract}

\section{Необходимые определения и полученные результаты}

\subsection{Двумерное пространство анти-де Ситтера}

Рассмотрим пространство  $\mathbb{R}^3_2 = \{ x = (x_1, x_2, x_3)\vert x_i \in \mathbb{R} \}$ с псевдоевклидовой метрикой $ds^2 = - dx_1^2-dx_2^2+dx_3^2$. Рассмотрим однополостный гиперболоид
$$
H_1^2 = \{x \in \mathbb{R}^3_2 \vert -x_1^2- x_2^2+x_3^2 = -1\},
$$
параметризуем его как
\begin{equation}
\label{par}
x_1 = \ch{\theta} \cos{\varphi}, \ x_2 = \ch{\theta} \sin{\varphi}, \ x_3 = \sh{\theta}, \ \theta \in \mathbb{R}, \ \varphi \in \mathbb{R}/(2 \pi \mathbb{Z}),
\end{equation}
и зададим на $H_1^2$ лоренцеву метрику $g = ds^2 \vert_{H_1^2}$.

\begin{definition}
{\it Двумерное пространство анти-де Ситтера} \cite{beem} есть односвязное накрывающее многообразие 
$$
\widetilde{H_1^2} = \{(\varphi, \theta) \in \mathbb{R}^2\},
$$
в которой задана лоренцева метрика $\tilde{g}$, индуцированная метрикой $g$.  
\end{definition}

Векторные поля $X_1, X_2 \in \VEC(\widetilde{H_1^2})$ образуют ортонормированный репер для метрики $\tilde{g}$, если
$$
\tilde{g}(X_2, X_2) = - \tilde{g}(X_1, X_1) = 1, \ \tilde{g}(X_1, X_2) = 0. 
$$

\subsection{Выражение метрики $\tilde{g}$ и ортонормированного репера  в координатах на $\widetilde{H_1^2}$}
Локально метрика $g$ совпадает с $\tilde{g}$. Метрика $g$ имеет такое выражение:
$$
g = d \theta^2 - \ch^2{\theta} d \varphi^2
$$

Ортонормированный репер:
$$
X_1 = \frac{1}{\ch{\theta}} \frac{\partial}{\partial \varphi},\ X_2 = \frac{\partial }{\partial \theta}
$$

\subsection{Постановка лоренцевой задачи на $\widetilde{H_1^2}$}
Лоренцевы длиннейшие для метрики $\tilde{g}$  суть решения следующей задачи оптимального управления:
\begin{center}
\begin{align}
& \dot q = u_1 X_1(q) + u_2 X_2(q), \qquad q = (\varphi, \theta) \in M = \widetilde{H_1^2}, \label{pr1}\\
& u \in U = \{(u_1, u_2) \in \mathbb{R}^2 \vert u_1^2 - u_2^2 \geq 0, \ u_1 > 0\}, \label{pr2}\\
& q(0) = q_0 = (\varphi, \theta_0), \qquad q(t_1) = q_1 = (\varphi, \theta_1), \label{pr3}\\
& J = \int_0^{t_1} \sqrt{u_1^2-u_2^2} \, dt \rightarrow \max. \label{pr4}
\end{align}
\end{center} 

\subsection{Условия принципа максимума Понтрягина для задачи на $\widetilde{H_1^2}$}

Гамильтониан ПМП, где $\nu \in \lbrace -1, 0 \rbrace$: 
$$
h_{u}^{\nu}(\lambda) = \langle \xi_1 d \theta+ \xi_2 d \varphi,  u_2 \frac{\partial}{\partial \theta} + \frac{u_1}{\ch{\theta}} \frac{\partial}{\partial \varphi} \rangle - \nu \sqrt{u_1^2-u_2^2} =
$$
$$
= \xi_1 u_2 + \xi_2 \frac{u_1}{\ch{\theta}} - \nu \sqrt{u_1^2-u_2^2} = h_1 u_1+h_2 u_2 - \nu \sqrt{u_1^2-u_2^2},
$$

где $h_1(\lambda) = \langle \lambda, X_1(q)\rangle = \frac{\xi_2}{\ch \theta}$, $h_2(\lambda) = \langle \lambda, X_2(q)\rangle = {\xi_1}$.

Условия принципа максимума Понтрягина:

\begin{enumerate}
\item Гамильтонова система

$$
\begin{cases}
\dot{\xi}_1 = -\partial_{\theta} h_{u}^{\nu}\\
\dot{\xi}_2 = -\partial_{\varphi} h_{u}^{\nu}\\
\dot{\theta} = \partial_{\xi_1} h_{u}^{\nu}\\
\dot{\varphi} = \partial_{\xi_2} h_{u}^{\nu}
\end{cases}
\Leftrightarrow
\begin{cases}
\dot{\xi}_1 = \xi_2 u_1 \frac{\sh{\theta}}{\ch^2{\theta}}\\
\dot{\xi}_2 = 0\\
\dot{\theta} = u_2\\
\dot{\varphi} = \frac{u_1}{\ch{\theta}}
\end{cases}
$$

\item Условие максимальности
$$
\xi_1 u_2 + \xi_2 \frac{u_1}{\ch{\theta}} - \nu \sqrt{u_1^2-u_2^2} \rightarrow \underset{u \in U}{\max}
$$

\item Нетривиальность
$$
(\xi_1,\xi_2, \nu) \neq (0,0,0)
$$
\end{enumerate}

\subsection{Анормальные траектории}

Анормальные траектории - траектории при $\nu = 0$.\\

Получаем такие случаи, когда максимум гамильтониана существует:

\begin{enumerate}
\item $h_1 = h_2 = h < 0$\\
$
\underset{u \in U}{\max} h_u^0 = 0\ \text{при}\ u_1 = -u_2
$

\item $h_1 = -h_2 = h < 0$\\
$
\underset{u \in U}{\max} h_u^0 = 0\ \text{при}\ u_1 = u_2
$
\end{enumerate}

Будем считать, что $u_1 = 1$, так как не следим за параметризацией кривой.

\begin{enumerate}
\item $h_1 = h_2 = h < 0,\ u_1 = 1 = -u_2$

$$
\varphi (t) = -\arctg{ \lbrace \sh{ (-t + \theta_0) } \rbrace} + \varphi_0 + \arctg{\lbrace \sh{\theta_0} \rbrace},\ \theta (t) =-t + \theta_0
$$

\item $h_1 = -h_2 = h < 0,\ u_1 = 1 = u_2$

$$
\varphi (t) = \arctg{ \lbrace \sh{( t + \theta_0 )} \rbrace} + \varphi_0 - \arctg{\lbrace \sh{\theta_0} \rbrace},\ \theta (t) =t + \theta_0
$$

\end{enumerate}

\subsection{Нормальные траектории}

Нормальные траектории - траектории при $\nu = -1$.

Получаем такие случаи, когда максимум гамильтониана существует:
$$
h_1 = -\ch{\beta},\ h_2 = \sh{\beta},\ u_1 = \rho \ch{\beta},\ u_2 = \rho \sh{\beta},\ \rho > 0
$$

Выберем на этой кривой параметризацию таким образом, что
$$
u_1^2 - u_2^2 = 1,\ u_1 > 0 \Rightarrow u_1 = \sqrt{u_2^2 + 1},\ h_1 = -\sqrt{u_2^2+1},\ h_2 = u_2
$$

Решение гамильтоновой системы выглядит так.\\

Есть первый интеграл 
$$
\ch{\psi} \ch{\theta} = D \geqslant 1
$$

Случай $D=1.$\\
Так как $\ch{x}\geqslant 1$, то оба гиперболических косинуса равны 1, поэтому $\theta \equiv 0, \psi \equiv 0$.
Тогда $\varphi \equiv -C_2t + D_1$.\\

Случай $D>1.$

$$
\varphi = \arctg{\Big( D \tg{(t+C_0)} \Big)} +E_0,\ \theta = \arsh{\Big( \pm \sqrt{D^2 - 1}\sin{(t+C_0) \Big)}},\ \psi(t) = \arsh{\Big( \pm \frac{\sqrt{D^2-1} \cos{(t+C_0)}}{\sqrt{1+(D^2-1)\sin^2{(t+C_0)}}} \Big)}
$$

Решение при $C_0 = 0,\ \theta(0) = 0,\ \varphi(0) = 0$: 
\begin{align*}
&\psi(t) = \arsh{\Big( \frac{s_0 \cos{t}}{\sqrt{1+s_0^2\sin^2{t}}} \Big)} \\
&\theta(t) = \arsh{\Big( s_0 \sin{t} \Big)} \\
&\varphi(t) = \arctg{\Big( c_0 \tg{t} \Big)}, 
\end{align*}
где $c_0=\ch{\psi(0)},\ s_0 = \sh{\psi(0)}$.

\subsection{Продолжение решения при $t>\pi/2$}

Посмотрев на полученную формулу для $\varphi$, можно подумать, что оно стремится к бесконечности за конечное время, но правая часть - гладкая ограниченная функция для всех $t \geqslant 0$, следовательно, и решение - гладкая функция для всех $t \geqslant 0$. Мы получаем такие формулы для продолжения:

$$
\begin{cases}
\theta(t) = \arsh{\Big( \sh{(\psi(0))}\sin{t} \Big)}\\
\varphi(t) = \begin{cases}
n\pi + \varphi_1(t - n\pi),\ t \in [n\pi, (2n+1)\pi/2]\\
(n+1)\pi - \varphi_1((n+1)\pi-t),\ t \in [(2n+1)\pi/2, (n+1)\pi]
\end{cases}
\end{cases}
$$

\begin{figure}[H]
\label{pic4}
\caption{График решения при $\psi(0) = 1$ }
\includegraphics[width = 10 cm, height = 8 cm]{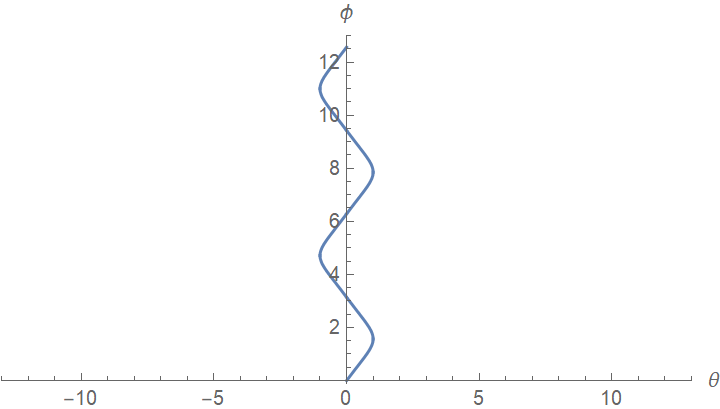}
\centering
\end{figure}

\begin{figure}[H]
\label{pic5}
\caption{Графики решений при $\psi(0) \in \lbrace 0,\ \pm 2,\ \pm 4,\ \pm 6,\ \pm 8,\ \pm10 \rbrace$}
\includegraphics[width = 10 cm, height = 8 cm]{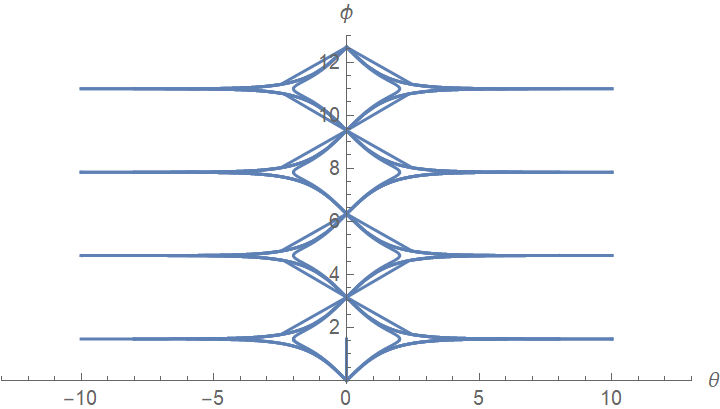}
\centering
\end{figure}

\subsection{Экспоненциальное отображение и его свойства}

Экспоненциальное отображение, определяемое следующими формулами: 
$$
\begin{cases}
\theta(t) = \arsh{\Big( \sh{(\psi(0))}\sin{t} \Big)}\\
\varphi(t) = \begin{cases}
\arctg{\Big( \ch{(\psi(0))}\tg{t} \Big)},\ t \in (0, \pi/2)\\
\pi/2,\ t = \pi/2 \\
\pi - \arctg{\Big( \ch{(\psi(0))}[\tg{(\pi - t)}] \Big)},\ t \in (\pi/2, \pi),
\end{cases}
\end{cases}
$$
задаёт диффеоморфизм областей\\
$A' = \lbrace (\psi_0,\ t_0) : \psi_0 \in \mathbb{R},\  t_0 \in (0,\pi) \rbrace$ и $C' = \lbrace (\theta, \varphi) \in \mathbb{R}^2 : \arctg{\sh{\theta}} < \varphi < \pi - \arctg{\sh{\theta}} \rbrace $.

\subsection{Множество достижимости из точки $\theta = 0,\ \varphi = 0$}

Множество достижимости это:
$$
V = \lbrace (\theta, \varphi) \in \mathbb{R}^2 : \varphi \geqslant \arctg{(\sh{\vert \theta \vert})} \rbrace
$$

\subsection{Область существования оптимальных траекторий}

Множество достижимости оптимальных траекторий из произвольной $q_0$ это множество:
$$
\mathcal B_{q_0} = \{(\theta, \varphi) \in M \vert \varphi_0 + \pi - \arctg{ \sh{ |\theta-\theta_0|}} >  \varphi \geqslant \varphi_0 + \arctg{ \sh {|\theta-\theta_0|}}\}
$$

\subsection{Оптимальный синтез на нижней границе области $C'$}

Пусть $\varphi_1 = \arcsin{(\sh{\theta_1})},\ (\varphi_1,\ \theta_1) \neq (0,0)$. Двигаться можно только по границе. Любой точке на нижней границе соответствуют анормальные траектории:
Если $\theta_1 > 0$, то выбираем эту формулу:
$$
\theta(t) = t,\ \varphi(t) = \arcsin{\sh{(t)}},\ u_1 = u_2
$$
Если $\theta_1 < 0$, то выбираем эту формулу:
$$
\theta(t) = -t,\ \varphi(t) = -\arcsin{\sh{(-t)}},\ u_1 = -u_2
$$

Отметим, что их длина равна $0$, так как $u_1^2-u_2^2 = 0$.

\subsection{Оптимальный синтез в области $C'$}

Момент, в который мы достигаем точки $q_1 = (\theta_1,\varphi_1) \in C'$:
$$
t_{q_1} = 
\begin{cases}
\arcsin{\sqrt{\frac{\tg^2{\varphi_1} - \sh^2{\theta_1}}{1+\tg^2{\varphi_1}}}},\ \varphi_1 \in (0,\pi/2)\\
\pi/2,\ \varphi_1 = \pi/2 \\
\pi - \arcsin{\sqrt{\frac{\tg^2{\varphi_1} - \sh^2{\theta_1}}{1+\tg^2{\varphi_1}}}},\ \varphi_1 \in (\pi/2,\pi)
\end{cases}
$$

Оптимальная траектория, соединяющая точки $(0,0)$ и $q_1 = (\theta_1,\varphi_1) \in C'$.\\

Если $t_{q_1} < \pi/2$, то

$$
\begin{cases}
\theta(t) = \arsh{\Big( \sh{(\psi_{q_1})}\sin{t} \Big)}\\
\varphi(t) = 
\arctg{\Big( \ch{(\psi_{q_1})}\tg{t} \Big)},\ t \in (0, t_{q_1})
\end{cases}
$$

Если $t_{q_1} = \pi/2$, то

$$
\begin{cases}
\theta(t) = \arsh{\Big( \sh{(\psi_{q_1})}\sin{t} \Big)}\\
\varphi(t) = 
\begin{cases}
\arctg{\Big( \ch{(\psi_{q_1})}\tg{t} \Big)},\ t \in (0, t_{q_1})\\
\pi/2
\end{cases}
\end{cases}
$$

Если $\pi/2 < t_{q_1} < \pi$, то

$$
\begin{cases}
\theta(t) = \arsh{\Big( \sh{(\psi_{q_1})}\sin{t} \Big)}\\
\varphi(t) = \begin{cases}
\arctg{\Big( \ch{(\psi_{q_1})}\tg{t} \Big)},\ t \in (0, \pi/2)\\
\pi/2,\ t = \pi/2\\
\pi - \arctg{\Big( \ch{(\psi_{q_1})}[\tg{(\pi - t)}] \Big)},\ t \in (\pi/2, t_{q_1}),
\end{cases}
\end{cases}
$$
где 
$$
\psi_{q_1} = \arsh{\Bigg( \sh{\theta_1}\sqrt{\frac{1+\tg^2{\varphi_1}}{\tg^2{\varphi_1} - \sh^2{\theta_1}}}\Bigg)}
$$

Отметим, что лоренцево расстояние до точки $q_1$ равно $t_{q_1}$, так как $u_1^2-u_2^2 = 1$.

\subsection{Программа для Wolfram Mathematica для вычисления траектории в области $C'$}

Некоторые графики, полученные с помощью написанной программы. Код написан в секции 2.13.

\begin{figure}[H]
\caption{Графики оптимальных траекторий}
\centering
\begin{subfigure}[H]{0.4\textwidth}
\label{grpr1}
\caption{Пример 1}
\includegraphics[width = 6 cm, height = 4 cm]{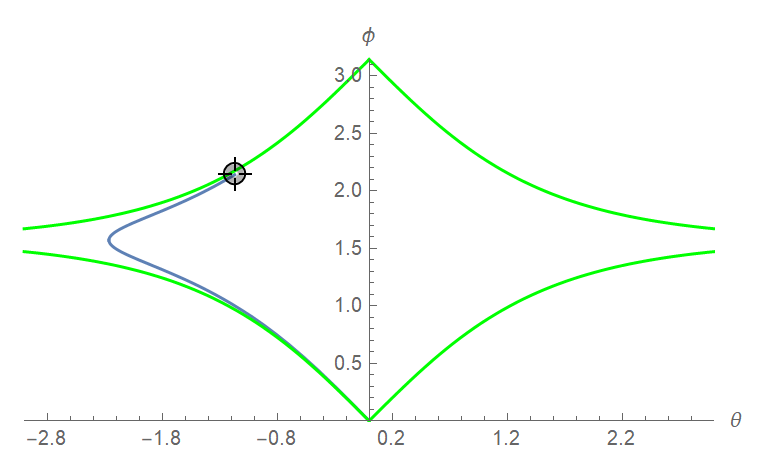}
\end{subfigure}

\begin{subfigure}[H]{0.4\textwidth}
\label{grpr2}
\caption{Пример 2}
\includegraphics[width = 6 cm, height = 4 cm]{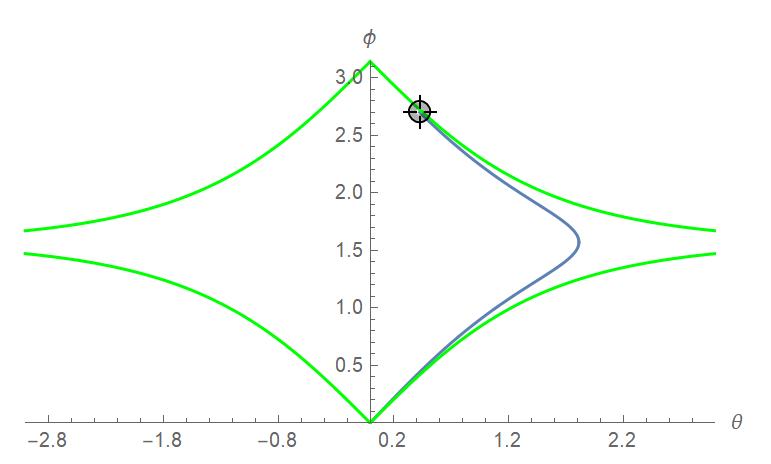}
\end{subfigure}

\begin{subfigure}[H]{0.4\textwidth}
\label{grpr3}
\caption{Пример 3}
\includegraphics[width = 6 cm, height = 4 cm]{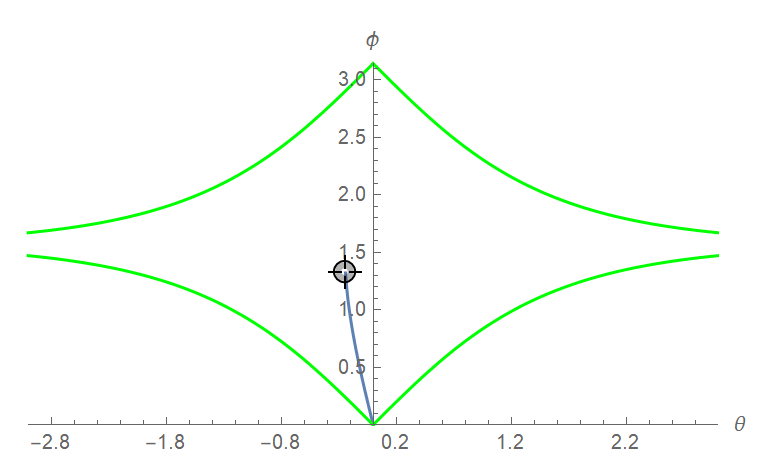}
\end{subfigure}
\end{figure}

\subsection{Точки выше верхней границы $C'$}

Для точек $(\theta_1, \varphi_1)$, таких что $\varphi_1 > \pi - \arctg{(\sh{|\theta_1|})}$, не существует оптимальной траектории. Лоренцево расстояние до них равно $+\infty$.

\subsection{Точки верхней границы $C'$}

Для точек $(\theta_1,\varphi_1)$, таких что $\varphi = \pi - \arctg{(\sh{|\theta_1|})}$ лоренцево расстояние равно $\pi$. Но только для точки $(0,\pi)$ существуют оптимальные траектории, причём их континуум. Для других точек этой кривой нет оптимальных траекторий.

\subsection{Поля Киллинга}

\begin{definition}

Векторное поле называется {\it полем Киллинга}, если производная Ли метрики вдоль него равно нулю.

\end{definition}

Векторное поле Киллинга является таковым тогда и только тогда, когда удовлетворяет уравнению: 

$$
X(g(V,W)) = g([X,V],W)+g(V,[X,W]),
$$

где $V$, $W$ - векторные поля, а $g$ - метрика \cite{lor_lob}.\\

Также векторные поля Киллинга образуют подалгебру Ли в алгебре Ли всех векторных полей на многообразии, причём для Лоренцева многообразия постоянной кривизны размерность равна $\frac{n(n+1)}{2}$, где $n$ - размерность многообразия \cite{lor_lob}.\\

В нашем случае алгебра полей Киллинга 3-мерна и её базисные вектора:
$$
\hat{X}_1 = \ch{\theta}X_1 = \partial_{\varphi},\ \hat{X}_2 = \sh{\theta}\cos{\varphi}X_1+\sin{\varphi}X_2,\ \hat{X}_3 = -\sh{\theta}\sin{\varphi}X_1+\cos{\varphi}X_2 
$$

Для них верны такие соотношения: 
$$
[\hat{X}_1,\hat{X}_2] = \hat{X}_3,\ [\hat{X}_2,\hat{X}_3]=-\hat{X}_1,\ [\hat{X}_3,\hat{X}_1]=\hat{X}_2
$$
Поэтому алгебра полей Киллинга изоморфна $\mathfrak{sl}(2)$ \cite{jakobson}. 

\subsection{Фазовые портреты полей Киллинга}

\begin{figure}[H]
\caption{Траектории полей Киллинга}
\centering
\begin{subfigure}[H]{0.4\textwidth}
\label{pic6}
\caption{Траектории $\hat{X}_1$}
\includegraphics[width = 6 cm, height = 6 cm]{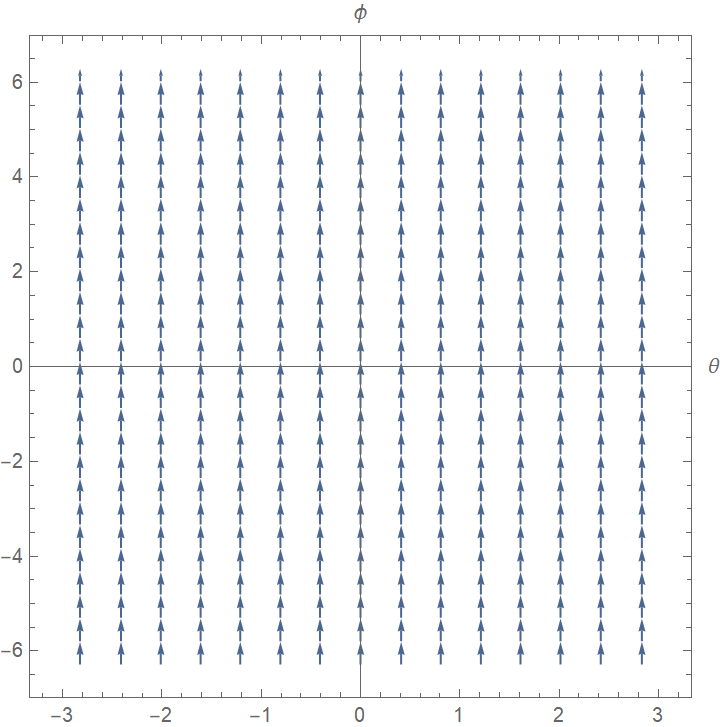}
\end{subfigure}

\begin{subfigure}[H]{0.4\textwidth}
\label{pic7}
\caption{Траектории $\hat{X}_2$}
\includegraphics[width = 6 cm, height = 6 cm]{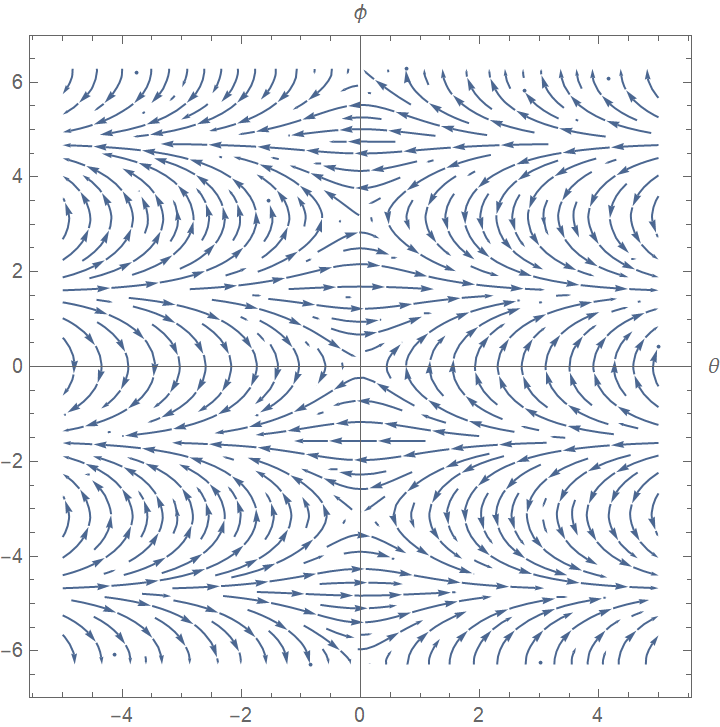}
\end{subfigure}

\begin{subfigure}[H]{0.4\textwidth}
\label{pic8}
\caption{Траектории $\hat{X}_3$}
\includegraphics[width = 6 cm, height = 6 cm]{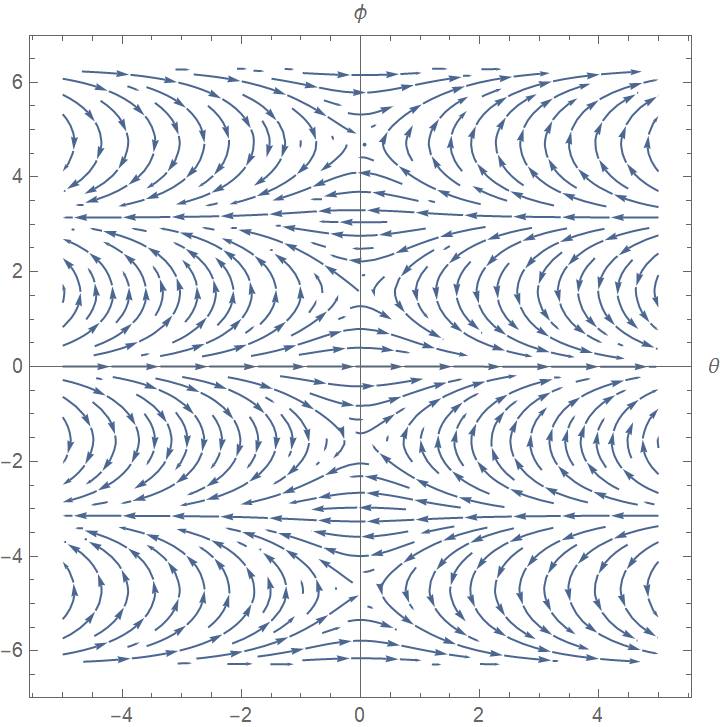}
\end{subfigure}
\end{figure}

\subsection{Как посчитать лоренцево расстояние между произвольными двумя точками $\widetilde{H_1^2}$?}

Предположим, что мы хотим посчитать расстояние между точками $q_0$ и $q_1$. Тогда мы перенесём точку $q_0$ в $(0,0)=q'_0$ вдоль траекторий полей Киллинга по определённому маршруту, а точку $q_1$ в точку $q'_1$ параллельно этому маршруту. Тогда $d(q_0,q_1) = d(q'_0,q'_1)$, так как значение метрики не меняется при переносе вдоль полей Киллинга. Мне удалось пока что только построить маршрут движения. Траектории полей Киллинга вычислить пока не получилось.

\section{Вычисления и обоснования}

В этой части номер подсекции совпадает с номером подсекции первой части, для обоснования результатов которой приведены соответствующие вычисления и соображения.

\subsection{Двумерное пространство анти-де Ситтера}

В этой подсекции нет фактов, которые надо обосновывать.

\subsection{Выражение метрики $\tilde{g}$ и ортонормированного репера  в координатах на $\widetilde{H_1^2}$}

Первым делом найдём $g$. Она локально совпадает с $\tilde{g}$.\\
$$
\dif{x_1} = \sh {\theta} \cos{\varphi} \dif{\theta} - \ch{\theta} \sin{\varphi} \dif{\varphi},\ \dif{x_2} = \sh{\theta} \sin{\varphi} \dif{\theta} + \ch{\theta} \cos{\varphi} \dif{\varphi},\ \dif{x_3} = \ch{\theta} \dif{\theta}
$$
Тогда 
$$
g = -\Big(\sh {\theta} \cos{\varphi} \dif{\theta} - \ch{\theta} \sin{\varphi} \dif{\varphi} \Big)^2 - \Big(\sh{\theta} \sin{\varphi} \dif{\theta} + \ch{\theta} \cos{\varphi} \dif{\varphi}\Big)^2 + \Big(\ch{\theta} \dif{\theta}\Big)^2=
$$
$$
= -\Big(\sh^2{\theta} \cos^2{\varphi} \dif{\theta}^2 - \sh {\theta} \cos{\varphi} \ch{\theta} \sin{\varphi}( \dif{\theta} \dif{\varphi} + \dif{\varphi}\dif{\theta} ) + \ch^2{\theta} \sin^2{\varphi} \dif{\varphi}^2 \Big)-
$$
$$
-\Big(\sh^2{\theta} \sin^2{\varphi} \dif{\theta}^2 + \sh{\theta} \sin{\varphi}\ch{\theta} \cos{\varphi}(\dif{\theta} \dif{\varphi} + \dif{\varphi}\dif{\theta} ) + \ch^2{\theta} \cos^2{\varphi} \dif{\varphi}^2 \Big)+
$$
$$
+\ch^2{\theta} \dif{\theta}^2 = -\sh^2{\theta}( \cos^2{\varphi} + \sin^2{\varphi} )\dif{\theta}^2 - \ch^2{\theta} ( \sin^2{\varphi} + \cos^2{\varphi} )\dif{\varphi}^2 + \ch^2{\theta} \dif{\theta}^2=
$$
$$
= (\ch^2{\theta} -\sh^2{\theta} )\dif{\theta}^2 - \ch^2{\theta}\dif{\varphi}^2 = \dif{\theta}^2 - \ch^2{\theta}\dif{\varphi}^2
$$

Теперь найдём собственные векторы, а потом откалибруем их длину по этой метрике.\\
Видно, что собственные числа $\lambda_1 = 1,$ $\lambda_2 = -\ch^2{\theta}$. Тогда собственные векторы, например, такие: 
$$
v_1 = \frac{\partial }{\partial \theta},\ v_2 = \frac{\partial}{\partial \varphi}
$$
Но
$$
g(v_1,v_1) = 1,\ g(v_2,v_2) = -\ch^2{\theta},\ g(v_1,v_2) = g(v_2,v_1) = 0,
$$
а нам надо:
$$
\tilde{g}(X_2, X_2) = - \tilde{g}(X_1, X_1) = 1, \qquad \tilde{g}(X_1, X_2) = 0
$$
Поэтому
$$
X_1 = \frac{1}{\ch{\theta}} \frac{\partial}{\partial \varphi},\ X_2 = \frac{\partial }{\partial \theta}
$$

\subsection{Постановка лоренцевой задачи на $\widetilde{H_1^2}$}

В этой подсекции нет фактов, которые надо обосновывать.

\subsection{Условия принципа максимума Понтрягина для задачи на $\widetilde{H_1^2}$}

В этой подсекции нет фактов, которые надо обосновывать.

\subsection{Анормальные траектории}

Анормальные траектории при $\nu = 0$:

Гамильтониан ПМП, где $\nu = 0$: 
$$
h_{u}^{0}(\lambda) = \xi_1 u_2 + \xi_2 \frac{u_1}{\ch{\theta}}
$$

Условия принципа максимума Понтрягина:

\begin{enumerate}
\item Гамильтонова система

$$
\begin{cases}
\dot{\xi}_1 = \xi_2 u_1 \frac{\sh{\theta}}{\ch^2{\theta}}\\
\dot{\xi}_2 = 0\\
\dot{\theta} = u_2\\
\dot{\varphi} = \frac{u_1}{\ch{\theta}}
\end{cases}
$$

\item Условие максимальности
$$
\xi_1 u_2 + \xi_2 \frac{u_1}{\ch{\theta}} \rightarrow \underset{u \in U}{\max}
$$

\item Нетривиальность
$$
(\xi_1,\xi_2) \neq (0,0)
$$
\end{enumerate}

Введ	ём линейные гамильтонианы $h_1(\lambda) = \langle \lambda, X_1(q)\rangle = \frac{\xi_2}{\ch{\theta}}$, 
$h_2(\lambda) = \langle \lambda, X_2(q)\rangle = {\xi_1}$.

\subsubsection{Условие максимальности}  

Записав условие максимальности через $h_1,\ h_2$: 
$$
h_1 u_1 + h_2 u_2 \rightarrow \underset{u \in U}{\max},
$$
поймём, когда максимум существует (принимает конечное значение) и при каких $h_1,\ h_2,\ u_1,\ u_2$ он достигается.\\
Рассмотрим такие случаи:

\begin{enumerate}

\item $h_2 \equiv 0,\ h_1 \neq 0$

$$
h_u^0 = h_1 u_1 \rightarrow 
\begin{cases}
+\infty,\ h_1 > 0,\ \text{при}\ u_1 \rightarrow +\infty\\
0,\ h_1 < 0,\ \text{при}\ u_1 \rightarrow 0
\end{cases}
$$
То есть максимум в данном случае не достигается.

\item $h_1 \equiv 0,\ h_2 \neq 0$

$$
h_u^0 = h_2 u_2 \rightarrow
\begin{cases}
+\infty,\ h_2 > 0,\ \text{при}\ u_2 \rightarrow +\infty\\
+\infty,\ h_2 < 0,\ \text{при}\ u_2 \rightarrow -\infty
\end{cases}
$$
То есть максимум в этом случае тоже не достигается.

\item $h_1 = h_2 = h \neq 0$
$$
h_u^0 = h(u_1 + u_2)
$$

$h>0$
$$
h_u^0 \rightarrow +\infty,\ u_1,\ u_2 \rightarrow +\infty
$$
То есть максимум в данном случае не достигается.\\
$h < 0$
$$
h_u^0 > 0\  \text{при}\ u_1 + u_2 < 0,\ \text{но так быть не может в силу определения множества}\ U
$$
Значит, что 
$$
\underset{u \in U}{\max} h_u^0 = 0\ \text{при}\ u_1 = -u_2
$$
\item $h_1 = -h_2 = h \neq 0$

$$
h_u^0 = h(u_1 - u_2)
$$

$h>0$

$$
h_u^0 \rightarrow +\infty \ \text{при}\ u_2 = 0,\ u_1 \rightarrow +\infty
$$
То есть максимум в данном случае не достигается.\\

$h<0$

$$
h_u^0 > 0\  \text{при}\ u_1 - u_2 < 0,\ \text{но так быть не может в силу определения множества}\ U
$$
Значит, что 
$$
\underset{u \in U}{\max} h_u^0 = 0\ \text{при}\ u_1 = u_2
$$

\item $h_1 > 0$
$$
h_u^0 = h_1 u_1 + h_2 u_2 \rightarrow +\infty \ \text{при}\ u_2 = 0,\ u_1 \rightarrow +\infty
$$
То есть максимум в данном случае не достигается.
\item $h_1 < 0,\ h_2 > 0,\ -h_1 < h_2$

$$
h_u^0 = h_1 u_1 + h_2 u_2 
$$
При $u_2 = u_1$
$$
h_u^0 = u_1(h_1 + h_2 ) \rightarrow +\infty \ \text{при}\ u_1 \rightarrow +\infty
$$
То есть максимум в данном случае не достигается.
\item $h_1 < 0,\ h_2 < 0,\ h_2 < h_1$

$$
h_u^0 = h_1 u_1 + h_2 u_2 \rightarrow +\infty \ \text{при}\ -u_2 = u_1 \rightarrow +\infty
$$
При $u_2 = -u_1$
$$
h_u^0 = u_1(h_1 - h_2 ) \rightarrow +\infty \ \text{при}\ u_1 \rightarrow +\infty
$$
То есть максимум в данном случае не достигается.
\item $h_1 < 0,\ h_2 > 0,\ h_2 < -h_1$

Так как 

$$
\begin{cases}
|h_1| > h_2 \\
u_1 \geqslant |u_2|
\end{cases}
\Rightarrow |h_1|u_1 \geqslant h_2 |u_2|,\ \text{но}\ h_1 < 0,\ h_2 > 0,\ \text{поэтому}
$$
максимум достигается при $u_2 = u_1$
$$
h_u^0 = u_1(h_1 + h_2) \rightarrow 0\ \text{при}\ u_1 \rightarrow 0
$$
То есть максимум в данном случае не достигается.
\item $h_1 <0,\ h_2 < 0,\ h_1 < h_2$
Так как 

$$
\begin{cases}
|h_1| > |h_2| \\
u_1 \geqslant |u_2|
\end{cases}
\Rightarrow |h_1|u_1 \geqslant |h_2| |u_2|,\ \text{но}\ h_1 < 0,\ h_2 < 0,\ \text{поэтому}
$$
максимум достигается при $u_2 = -u_1$
$$
h_u^0 = u_1(h_1 + h_2) \rightarrow 0\ \text{при}\ u_1 \rightarrow 0
$$
То есть максимум в данном случае не достигается.
\end{enumerate}
В итоге получаем такие случаи, когда максимум существует:

\begin{enumerate}
\item $h_1 = h_2 = h < 0$
$$
\underset{u \in U}{\max} h_u^0 = 0\ \text{при}\ u_1 = -u_2
$$

\item $h_1 = -h_2 = h < 0 $
$$
\underset{u \in U}{\max} h_u^0 = 0\ \text{при}\ u_1 = u_2
$$
\end{enumerate}

\subsubsection{Гамильтонова система}

Решим систему уравнений.

Из системы сразу следует, что 
$$
\xi_2 \equiv const = c_2
$$

Будем считать, что $u_1 = 1$, так как не следим за параметризацией кривой.

\begin{enumerate}

\item

$$
\begin{cases}
\dot{\xi}_1 = \xi_2 \frac{\sh{\theta}}{\ch^2{\theta}}\\
\dot{\xi}_2 = 0\\
\dot{\theta} = -1\\
\dot{\varphi} = \frac{1}{\ch{\theta}}
\end{cases}
$$

$$
h_2 = \xi_1 = \frac{\xi_2}{\ch{\theta}} = \frac{c_2}{\ch{\theta}} = h_1 < 0 \Rightarrow c_2 < 0,\ \xi_1 = \frac{c_2}{\ch{\theta}}
$$

$\xi_1,\ \xi_2$ мы нашли, теперь найдём $\theta,\ \varphi$. Разделим третье уравнение на четвёртое.

$$
\frac{\dif \theta}{\dif \varphi} = -\ch{\theta} \Leftrightarrow -\frac{\dif \theta}{\ch{\theta}} = \dif \varphi
$$

$$
\int \frac{\dif \theta}{\ch{\theta}} = \int \frac{\ch{\theta} \dif \theta}{\ch^2{\theta}}= \int \frac{\ch{\theta} \dif \theta}{1+\sh^2{\theta}} = \lbrace \sh{\theta} = \gamma,\ \ch{\theta} \dif \theta = \dif \gamma \rbrace = \int \frac{\dif \gamma}{\gamma^2 + 1} =
$$
$$
= \arctg{\gamma} + D_0
= \arctg{\lbrace \sh{\theta} \rbrace} + D_0
$$

Тогда
$$
-\arctg{\lbrace \sh{\theta} \rbrace} + D_0 = \varphi
$$

Так как $\theta(0) = \theta_0,\ \varphi (0) = \varphi_0$, то 
$$
D_0 = \varphi_0 + \arctg{\lbrace \sh{\theta_0} \rbrace}
$$
Так как $\dot{\theta} = -1$, то
$$
\theta (t) = \int_{0}^{t} (-1)\dif s + D_1 = -t + D_1
$$
Так как $\theta(0) = \theta_0$, то
$$
D_1 = \theta_0
$$
Поэтому
$$
\varphi (t) = -\arctg{ \lbrace \sh{ (-t + \theta_0) } \rbrace} + \varphi_0 + \arctg{\lbrace \sh{\theta_0} \rbrace},\ \theta (t) =-t + \theta_0
$$

\item

$$
\begin{cases}
\dot{\xi}_1 = \xi_2 \frac{\sh{\theta}}{\ch^2{\theta}}\\
\dot{\xi}_2 = 0\\
\dot{\theta} = 1\\
\dot{\varphi} = \frac{1}{\ch{\theta}}
\end{cases}
$$

$$
-\xi_1 = \frac{\xi_2}{\ch{\theta}} = \frac{c_2}{\ch{\theta}} < 0 \Rightarrow c_2 < 0,\ \xi_1 = -\frac{c_2}{\ch{\theta}}
$$

$\xi_1,\ \xi_2$ мы нашли, теперь найдём $\theta,\ \varphi$. Разделим третье уравнение на четвёртое.

$$
\frac{\dif \theta}{\dif \varphi} = \ch{\theta} \Leftrightarrow \frac{\dif \theta}{\ch{\theta}} = \dif \varphi
$$

Тогда
$$
\varphi = \arctg{\lbrace \sh{\theta} \rbrace } + D_2 
$$
Так как $\theta(0) = \theta_0,\ \varphi (0) = \varphi_0$, то
$$
D_2 = \varphi_0 - \arctg{\lbrace \sh{\theta_0} \rbrace }
$$
Так как $\dot{\theta} = 1$, то
$$
\theta (t) = \int_{0}^{t} \dif s + D_3 = t + D_3
$$
Так как $\theta(0) = \theta_0$, то
$$
D_3 = \theta_0
$$

Поэтому
$$
\varphi (t) = \arctg{ \lbrace \sh{( t + \theta_0 )} \rbrace} + \varphi_0 - \arctg{\lbrace \sh{\theta_0} \rbrace},\ \theta (t) =t + \theta_0
$$

\end{enumerate}

\subsection{Нормальные траектории}

Нормальные траектории при $\nu = -1$:

Гамильтониан ПМП, где $\nu = -1$: 
$$
h_{u}^{0}(\lambda) = \xi_1 u_2 + \xi_2 \frac{u_1}{\ch{\theta}} + \sqrt{u_1^2-u_2^2}
$$

Условия принципа максимума Понтрягина:

\begin{enumerate}
\item Гамильтонова система

$$
\begin{cases}
\dot{\xi}_1 = \xi_2 u_1 \frac{\sh{\theta}}{\ch^2{\theta}}\\
\dot{\xi}_2 = 0\\
\dot{\theta} = u_2\\
\dot{\varphi} = \frac{u_1}{\ch{\theta}}
\end{cases}
$$

\item Условие максимальности
$$
\xi_1 u_2 + \xi_2 \frac{u_1}{\ch{\theta}} + \sqrt{u_1^2-u_2^2} \rightarrow \underset{u \in U}{\max}
$$

\item Нетривиальность
$$
(\xi_1,\xi_2) \neq (0,0)
$$
\end{enumerate}

\subsubsection{Условие максимальности}

Записав условие максимальности через $h_1 = \frac{\xi_2}{\ch{\theta}},\ h_2 = \xi_1$: 
$$
h_1 u_1 + h_2 u_2 + \sqrt{u_1^2-u_2^2} \rightarrow \underset{u \in U}{\max},
$$
поймём, когда максимум существует (принимает конечное значение) и при каких $h_1,\ h_2,\ u_1,\ u_2$ он достигается.\\
Рассмотрим такие случаи:

\begin{enumerate}
\item $h_1 > 0,\ h_1 > |h_2|$

$$
\begin{cases}
u_1 = \rho \ch{\alpha}\\
u_2 = \rho \sh{\alpha}
\end{cases}
\text{и}\ 
\begin{cases}
h_1 = r\ch{\beta}\\
h_2 = r\sh{\beta}
\end{cases}
$$

$$
h_{u}^{-1} = r\rho \ch{\alpha}\ch{\beta} + r\rho \sh{\alpha}\sh{\beta} + \rho = \rho (r \ch{(\alpha+\beta)} + 1) \rightarrow +\infty \ \text{при}\ \rho \rightarrow +\infty
$$

В этом случае максимум не достигается.

\item $h_1 > 0,\ h_2 = \pm h_1 = \pm h$

$$
\begin{cases}
u_1 = \rho \ch{\alpha}\\
u_2 = \rho \sh{\alpha}
\end{cases}
$$

$$
h_{u}^{-1} = h \rho( \ch{\alpha} \pm \sh{\alpha} ) + \rho = \rho( h( \ch{\alpha} \pm \sh{\alpha} ) + 1 ) \rightarrow +\infty \ \text{при}\ \alpha = 0,\ \rho \rightarrow +\infty
$$

В этом случае максимум не достигается.

\item $|h_2| > h_1 \Leftrightarrow h_1^2 - h_2^2 < 0$

Введём координаты 
$$
\begin{cases}
u_1 = \rho \ch{\alpha}\\
u_2 = \rho \sh{\alpha}
\end{cases}
\text{и}\ 
\begin{cases}
h_1 = r\sh{\beta}\\
h_2 = \pm r\ch{\beta}
\end{cases}
$$

Тогда гамильтониан перепишется:
$$
h_{u}^{-1} = h_1 u_1 + h_2 u_2 + \sqrt{u_1^2 - u_2^2} = r\rho \sh{\beta}\ch{\alpha} \pm r\rho \sh{\alpha}\ch{\beta} + \rho = 
$$
$$
=
\begin{cases}
\rho(1+r\sh{(\alpha + \beta)})\\
\rho(1-r\sh{(\alpha-\beta)})
\end{cases}
$$

Заметим, что при, соответственно, $\alpha = -\beta,\ \alpha = \beta$ получаем, что $h_{u}^{-1} = \rho \Rightarrow$
$$
\underset{u \in U}{\max} h_{u}^{-1} \rightarrow +\infty \ \text{при}\ \rho \rightarrow +\infty
$$
В этом случае максимум не достигается.

\item $h_1 < 0,\ h_2 > 0,\ -h_2 = h_1 = h$

Запишем в координатах:

$$
\begin{cases}
u_1 = \rho \ch{\alpha}\\
u_2 = \rho \sh{\alpha}
\end{cases}
$$

$$
h_{u}^{-1} = h(u_1-u_2)+\sqrt{u_1^2-u_2^2} = h\rho(\ch{\alpha} - \sh{\alpha}) + \rho = \rho \Big(1+ \frac{h}{2}( e^{\alpha} + e^{-\alpha} - e^{\alpha} + e^{-\alpha} )\Big) = 
$$
$$
= \rho ( 1 + h e^{-\alpha} ) 
$$

Можем найти такой $\alpha$, что 
$$
e^{-\alpha} = \frac{1}{2|h|} \Leftrightarrow \alpha = -\ln{\frac{1}{2|h|}} \Leftrightarrow \alpha = \ln{2|h|}
$$
Тогда
$$
h_{u}^{-1} = \rho(1+h\frac{1}{2|h|}) = \frac{\rho}{2} \rightarrow +\infty \ \text{при}\ \rho \rightarrow +\infty
$$
Максимум в этом случае не достигается.

\item $h_1 < 0,\ h_2 < 0,\ h_2 = h_1=h$
Запишем в координатах:

$$
\begin{cases}
u_1 = \rho \ch{\alpha}\\
u_2 = \rho \sh{\alpha}
\end{cases}
$$

$$
h_{u}^{-1} = h(u_1+u_2) + \sqrt{u_1^2-u_2^2}=h\rho(\ch{\alpha} + \sh{\alpha}) + \rho = \rho \Big(1+ \frac{h}{2}( e^{\alpha} + e^{-\alpha} + e^{\alpha} - e^{-\alpha} )\Big) = 
$$
$$
= \rho ( 1 + h e^{\alpha} ) 
$$

Можем найти такой $\alpha$, что 
$$
e^{\alpha} = \frac{1}{2|h|} \Leftrightarrow \alpha = \ln{\frac{1}{2|h|}}
$$
Тогда
$$
h_{u}^{-1} = \rho(1+h\frac{1}{2|h|}) = \frac{\rho}{2} \rightarrow +\infty \ \text{при}\ \rho \rightarrow +\infty
$$
Максимум в этом случае не достигается.

\item $h_1 < 0,\ h_2 < |h_1| \Leftrightarrow h_1^2 - h_2^2 > 0$

Введём новые координаты:

$$
\begin{cases}
u_1 = \rho \ch{\alpha}\\
u_2 = \rho \sh{\alpha}
\end{cases}
\text{и}\ 
\begin{cases}
h_1 = -r\ch{\beta}\\
h_2 = r\sh{\beta}
\end{cases}
\text{где}\ \rho,\ r \geqslant 0,\ \alpha,\ \beta \in \mathbb{R}
$$

Тогда гамильтониан перепишется: 

$$
h_{u}^{-1} = h_1 u_1 + h_2 u_2 + \sqrt{u_1^2 - u_2^2} = -r \rho \ch{\beta} \ch{\alpha} + r \rho \sh{\beta} \sh{\alpha} + \sqrt{\rho^2(\ch^2{\alpha} - \sh^2{\alpha})} = 
$$
$$
= -r \rho( \ch{\beta} \ch{\alpha} - \sh{\beta} \sh{\alpha} ) + \rho = -r \rho \ch{(\alpha - \beta)} + \rho = \rho \Big( 1 - r \ch{(\alpha - \beta)} \Big)
$$

$\rho$ и $\alpha$ независимы, поэтому максимизировать произведение надо, максимизируя каждый из сомножителей.\\ Чтобы максимизировать скобку $1-r\ch{(\alpha - \beta)}$, надо минимизировать $\ch{(\alpha-\beta)}$. Функция $\ch{(\cdot)}$ принимает минимальное значение в нуле, поэтому получаем условие $\alpha = \beta$. Теперь скобка выглядит $1-r$.

Получаем 3 случая:
$$
1-r > 0 \Rightarrow \underset{u \in U}{\max} h_{u}^{-1} \rightarrow +\infty \ \text{при}\ \rho \rightarrow +\infty
$$
Максимум в этом случае не достигается.
$$
1-r < 0 \Rightarrow \underset{u \in U}{\max} h_{u}^{-1} \rightarrow 0 \ \text{при}\ \rho \rightarrow 0
$$
Максимум в этом случае не достигается.
$$
1-r = 0 \Rightarrow \underset{u \in U}{\max} h_{u}^{-1} = 0
$$
Максимум в этом случае достигается.

\end{enumerate}

\subsubsection{Решение гамильтоновой системы}

В итоге получаем такие случаи, когда максимум существует:
$$
h_1 = -\ch{\beta},\ h_2 = \sh{\beta},\ u_1 = \rho \ch{\beta},\ u_2 = \rho \sh{\beta},\ \rho > 0
$$

Иначе
$$
h_1^2 - h_2^2 = 1,\ h_1 < 0,\ u_1 = -\rho h_1,\ u_2 = \rho h_2,\ \rho > 0 \Rightarrow u_1^2 - u_2^2 = \rho^2;\ h_1 = -\frac{u_1}{\rho},\ h_2 = \frac{u_2}{\rho}
$$

Выберем на этой кривой параметризацию таким образом, что
$$
u_1^2 - u_2^2 = 1,\ u_1 > 0 \Rightarrow u_1 = \sqrt{u_2^2 + 1},\ h_1 = -\sqrt{u_2^2+1},\ h_2 = u_2
$$

$$
h_1 = \frac{\xi_2}{\ch{\theta}} = -\sqrt{u_2^2+1},\ h_2 = \xi_1 = u_2 \Rightarrow \xi_2 = -\ch{\theta}\sqrt{u_2^2+1},\ \xi_1 = u_2 \Rightarrow
$$
$$
\Rightarrow \dot{\xi}_1 = \dot{h}_2,\ h_1 \ch{\theta} = \xi_2,\ u_1 = -h_1,\ u_2 = h_2
$$

$$
\begin{cases}
\dot{\xi}_1 = \xi_2 u_1 \frac{\sh{\theta}}{\ch^2{\theta}}\\
\dot{\xi}_2 = 0\\
\dot{\theta} = u_2\\
\dot{\varphi} = \frac{u_1}{\ch{\theta}}
\end{cases}
$$

Из второго уравнения имеем $\xi_2 = C_2 \equiv const \Rightarrow h_1 = \frac{C_2}{\ch{\theta}}$.\\

А также запишем $h_1 = -\ch{\psi},\ h_2 = \sh{\psi} \Rightarrow \dot{h}_2 = \dot{\psi} \ch{\psi}$\\
Поэтому система записывается так:
$$
\begin{cases}
\dot{\psi} \ch{\psi} = -C_2h_1\frac{\sh{\theta}}{\ch^2{\theta}} = -h_1^2 \frac{\sh{\theta}}{\ch{\theta}} = -\ch^2{\psi} \frac{\sh{\theta}}{\ch{\theta}}\\
\dot{\theta} = \sh{\psi}\\
\dot{\varphi} = -\frac{C_2}{\ch^2{\theta}} = \frac{\ch{\psi}}{\ch{\theta}}
\end{cases}
\Rightarrow
$$

$$
\Rightarrow 
\dot{\psi} = -\ch{\psi} \frac{\sh{\theta}}{\ch{\theta}},\ \dot{\theta} = \sh{\psi}
$$
Поделим первое уравнение на второе:
$$
\frac{\dif \psi}{\dif \theta} = -\cth{\psi} \th{\theta} \Leftrightarrow \th{\psi} \dif{\psi} = -\th{\theta} \dif{\theta} \Leftrightarrow \ln{(\ch{\psi})} = -\ln{(\ch{\theta})} + \ln{D} \Leftrightarrow
$$
$$
\Leftrightarrow \ch{\psi}\ch{\theta} \equiv const = D \geqslant 1 
$$

Случай $D=1.$\\
Так как $\ch{x}\geqslant 1$, то оба гиперболических косинуса равны 1, поэтому $\theta \equiv 0, \psi \equiv 0$.
Тогда $\varphi \equiv -C_2t + D_1$.\\

Случай $D>1.$
$$
\ch{\psi}\ch{\theta} = D \Rightarrow \ch{\psi} = \frac{D}{\ch{\theta}} \Rightarrow 
\sh{\psi} =  \pm \sqrt{\frac{D^2}{\ch^2{\theta}} - 1}
$$

$$
\dot{\theta} = \pm \sqrt{\frac{D^2}{\ch^2{\theta}} - 1} \Leftrightarrow \frac{\ch{\theta} \dif \theta}{\sqrt{D^2 - \ch^2{\theta}}} = \pm dt
$$
$$
\int \frac{\ch{\theta} \dif \theta}{\sqrt{D^2 - \ch^2{\theta}}} = \lbrace \sh{\theta} = \alpha,\ \ch{\theta} \dif \theta  = \dif \alpha \rbrace = \int \frac{\dif \alpha}{\sqrt{D^2 - 1-\alpha^2}} = \arcsin{\frac{\sh{\theta}}{\sqrt{D^2 - 1}}} - C_0
$$

$$
\arcsin{\frac{\sh{\theta}}{\sqrt{D^2 - 1}}} = \pm(t+C_0) \Leftrightarrow \frac{\sh{\theta}}{\sqrt{D^2 - 1}} = \pm \sin{(t+C_0)} \Leftrightarrow \sh{\theta} = \pm \sqrt{D^2 - 1}\sin{(t+C_0)} \Leftrightarrow
$$
$$
\theta = \arsh{\Big( \pm \sqrt{D^2 - 1}\sin{(t+C_0) \Big)}}
$$
Тем самым нашли функцию $\theta (t)$. \\
Найдём теперь функцию $\psi(t)$:
$$
\sh{\psi} = \pm \sqrt{\frac{D^2}{\ch^2{\theta}} - 1} = \pm \sqrt{\frac{D^2 - 1-\sh^2{\theta}}{1+\sh^2{\theta}}} = \pm \sqrt{\frac{D^2 - 1 - (D^2-1)\sin^2{(t+C_0)}}{1+(D^2-1)\sin^2{(t+C_0)}}} =
$$
$$
= \pm \sqrt{\frac{D^2-1 - (D^2-1) + (D^2-1)\cos^2{(t+C_0)}}{1+(D^2-1)\sin^2{(t+C_0)}}} = \pm \frac{\sqrt{D^2-1} \cos{(t+C_0)}}{\sqrt{1+(D^2-1)\sin^2{(t+C_0)}}}
$$
Тогда
$$
\psi(t) = \arsh{\Big(\pm  \frac{\sqrt{D^2-1} \cos{(t+C_0)}}{\sqrt{1+(D^2-1)\sin^2{(t+C_0)}}} \Big)}
$$
Осталось найти $\varphi(t)$:
\begin{equation}
\label{dphi}
\dot{\varphi} = -\frac{C_2}{\ch^2{\theta}} = \frac{\ch{\psi}}{\ch{\theta}} = \frac{D}{\ch^2{\theta}} = \frac{D}{1+\sh^2{\theta}} = \frac{D}{1+(D^2-1)\sin^2{(t+C_0)}} \Leftrightarrow
\end{equation}
$$
\Leftrightarrow \dif \varphi = \frac{D \dif t}{1+(D^2-1)\sin^2{(t+C_0)}}
$$
Найдём соответствующий правой части интеграл:
$$
\int \frac{\dif t}{a^2\sin^2{t} + 1} = \lbrace \tg{t} = \alpha,\ \frac{\dif t}{\cos^2{t}} = \dif \alpha,\ \frac{1}{\cos^2{t}} = 1 + \tg^2{t},\ \sin^2{t} = \frac{\tg^2{t}}{1+\tg^2{t}} \rbrace =
$$
$$
= \int \frac{\dif \alpha}{(1+\alpha^2)(\frac{a^2 \alpha^2}{1+\alpha^2} + 1)} = \int \frac{\dif \alpha}{(a^2 +1) \alpha^2 +1} = \frac{1}{\sqrt{a^2 + 1}}\arctg{( \alpha \sqrt{a^2+1} )} + E_0
$$
Отсюда получаем, что
$$
\varphi = \frac{D}{\sqrt{D^2}}\arctg{\Big( \sqrt{D^2}} \tg{(t+C_0)  \Big)} + E_0 = \arctg{\Big( D \tg{(t+C_0)} \Big)} +E_0
$$
$$
\theta = \arsh{\Big( \pm \sqrt{D^2 - 1}\sin{(t+C_0) \Big)}},\ \psi(t) = \arsh{\Big( \pm \frac{\sqrt{D^2-1} \cos{(t+C_0)}}{\sqrt{1+(D^2-1)\sin^2{(t+C_0)}}} \Big)}
$$

\subsubsection{Решение с начальными условиями $\theta(0) = \theta_0$, $\varphi(0)=0$}
Начальные условия:
$$
\theta(0) =0,\ \varphi(0) = 0
$$

$$
0 = \varphi{(0)} = E_0
$$

$$
\ch{\psi{(0)}} = \frac{D}{\ch{\theta{(0)}}} = D,\ \text{а}\ \sqrt{D^2-1} = \pm \sh{\psi{(0)}}\ \text{в зависимости от знака}\ \psi{(0)}
$$

$$
0 = \theta{(0)}= \arsh{(\pm \sqrt{D^2-1}\sin{C_0})} \Rightarrow C_0 \in \lbrace k\pi \vert k \in \mathbb{Z} \rbrace \ \text{или}\ \psi(0) = 0
$$

Рассмотрим набор функций при $C_0 = 0$, $\psi{(0)} \geqslant 0$:
$$
\begin{cases}
\psi(t) = \arsh{\Big( \frac{\sh{(\psi(0))} \cos{t}}{\sqrt{1+\sh^2{(\psi(0))}\sin^2{t}}} \Big)}\\
\theta(t) = \arsh{\Big( \sh{(\psi(0))}\sin{t} \Big)}\\
\varphi(t) = \arctg{\Big( \ch{(\psi(0))}\tg{t} \Big)}
\end{cases}
$$

Этот набор функций является решением системы дифференциальных уравнений с начальными условиями $\theta(0)=0,\ \varphi(0) = 0,\ \psi(0) \geqslant 0$, что следует из сделанных выше выкладок:
$$
\begin{cases}
\dot{\psi} =  -\ch{\psi} \frac{\sh{\theta}}{\ch{\theta}}\\
\dot{\theta} = \sh{\psi}\\
\dot{\varphi} = \frac{\ch{\psi}}{\ch{\theta}}
\end{cases}
$$

Рассмотрим систему аналитических дифференциальных уравнений, где $z \in \mathbb{C}$:
$$
\begin{cases}
\psi'(z) = -\ch{\psi(z)} \frac{\sh{\theta(z)}}{\ch{\theta(z)}}\\
\theta'(z) = \sh{\psi(z)}\\
\varphi'(z) = \frac{\ch{\psi(z)}}{\ch{\theta(z)}}
\end{cases}
$$

Буквальной проверкой можно увидеть, что набор функций, определённый выше, с заменой вещественного $t$ на комплексный $z$ становится решениями данной системы аналитических дифференциальных уравнений с тем же начальным условием. Пусть ряд Тейлора этого решения сходится в некотором поликруге $U_0$ с центром в точке $(0,0,\psi(0))$. Теперь построим путь, соединяющий неотрицательное $\psi(0)$ с отрицательным числом $\gamma_1$, осуществив аналитическое продолжение по цепочке поликругов (последний из них это поликруг $U_1$ с центром $(0,0,\gamma_1)$), причём наш путь должен быть таким, что правые части системы тоже вдоль него аналитически продолжаются. Тогда, по теореме 6.5.1 из [\ref{Bibik}, стр. 196], мы получим решение с заданным начальным условием ($\theta(0) = 0,\ \varphi(0) = 0,\ \psi(0)=\gamma_1 < 0$) в поликруге $U_1$. И так можно сделать для всех вещественных $\psi(0)$, поэтому решение исходного вещественного дифференциального уравнения выглядит так:\\
$$
\begin{cases}
\psi(t) = \arsh{\Big( \frac{\sh{(\psi(0))} \cos{t}}{\sqrt{1+\sh^2{(\psi(0))}\sin^2{t}}} \Big)}\\
\theta(t) = \arsh{\Big( \sh{(\psi(0))}\sin{t} \Big)}\\
\varphi(t) = \arctg{\Big( \ch{(\psi(0))}\tg{t} \Big)}
\end{cases}
$$

Обозначим $c_0=\ch{\psi(0)},\ s_0 = \sh{\psi(0)}$:
\begin{align}
&\psi(t) = \arsh{\Big( \frac{s_0 \cos{t}}{\sqrt{1+s_0^2\sin^2{t}}} \Big)} \nonumber \\
&\theta(t) = \arsh{\Big( s_0 \sin{t} \Big)}  \label{theta}\\
&\varphi(t) = \arctg{\Big( c_0 \tg{t} \Big)}. \label{phi} 
\end{align}
Обозначим также

\begin{align}
&X(t) = s_0 \sin{t} \label{X}\\
&Y(t) = c_0 \tg{t} \label{Y}
\end{align}

\subsection{Продолжение решения при $t>\pi/2$}

Вернёмся к уравнению гамильтоновой системы 
$$
\dot{\varphi} = \frac{D}{1+(D^2-1)\sin^2{t}},\ D > 1,\ \varphi(0) = 0
$$
Правая часть - гладкая ограниченная функция для всех $t \geqslant 0$, следовательно, и решение - гладкая функция для всех $t \geqslant 0$:
$$
\varphi(t) = \int_0^{t} \frac{D}{1+(D^2-1)\sin^2{\tau}} \,d \tau \
$$
Мы его проинтегрировали и получили формулу, верную на полуинтервале $t \in [0, \pi/2)$:
$$
\varphi_0(t) =  \arctg{\Big( D \tg{t} \Big)}
$$

С другой стороны, заметим, что производная является периодической функцией с периодом $\pi$.\\
В точках $k \pi,\ k \in \mathbb{Z}$ функция достигает максимального значения, а в точке $n \pi/2,\ n=2l+1, l \in \mathbb{Z}$ - минимального. Причём она чётная на каждом промежутке $[(k-1)\pi,k\pi],\ k \in \mathbb{Z}$ относительно середины $(2k-1)\pi/2$.\\
В силу непрерывности решения, попробуем склеить её на каждом промежутке
$$
\underset{t \rightarrow \pi/2-0}{\lim} \arctg{\Big( D \tg{t} \Big)} = \pi/2
$$
По непрерывности, $\underset{t \rightarrow \pi/2+0}{\lim} \varphi(t) = \pi/2$.\\
В силу чётности производной относительно точки $\pi/2$, продлим решение на отрезок $[\pi/2,\pi]$ по формуле:
$$
t \in [0,\pi/2],\ \varphi(\pi/2 + t) = \pi/2 + \pi/2 - \arctg{\Big( D \tg{(\pi/2 - t)} \Big)},
$$
где первый член из условия, что $\varphi(t) = \pi/2$, а $\pi/2 - \arctg{\Big( D \tg{(\pi/2 - t)} \Big)}=\pi/2 - \varphi_0 (\pi/2 - t) = \int_{\pi/2}^{\pi/2+t} \frac{D}{1+(D^2-1)\sin^2{\tau}} \, d \tau$.
Соответственно, $\varphi(\pi) = \pi$, и мы получаем общую формулу, где $n \in \mathbb{N} \cup 0$:
$$
\varphi(\tau) = 
\begin{cases}
n\pi + \varphi_0(\tau - n\pi),\ \tau \in [n\pi, (2n+1)\pi/2]\\
(n+1)\pi - \varphi_0((n+1)\pi-\tau),\ \tau \in [(2n+1)\pi/2, (n+1)\pi]
\end{cases}
$$
Вернёмся к формулам для решения:
$$
\begin{cases}
\theta(t) = \arsh{\Big( \sh{(\psi(0))}\sin{t} \Big)}\\
\varphi(t) = \arctg{\Big( \ch{(\psi(0))}\tg{t} \Big)}
\end{cases}
$$
и запишем новую функцию $\varphi(t),\ t \in [0,+\infty)$ через $\varphi_1(t) =  \arctg{\Big( \ch{(\psi(0))}\tg{t} \Big)}$:

\begin{equation}
\begin{cases}
\theta(t) = \arsh{\Big( \sh{(\psi(0))}\sin{t} \Big)}\\
\varphi(t) = \begin{cases}
n\pi + \varphi_1(t - n\pi),\ t \in [n\pi, (2n+1)\pi/2]\\ \label{full_resh}
(n+1)\pi - \varphi_1((n+1)\pi-t),\ t \in [(2n+1)\pi/2, (n+1)\pi]
\end{cases}
\end{cases}
\end{equation}

\subsection{Экспоненциальное отображение и его свойства}

Покажем, что $A' = \lbrace (\psi_0,\ t_0) : \psi_0 \in \mathbb{R},\  t_0 \in (0,\pi) \rbrace$ диффеоморфно $C' = \lbrace (\theta, \varphi) \in \mathbb{R}^2 : \arctg{\sh{\theta}} < \varphi < \pi - \arctg{\sh{\theta}} \rbrace $.\\
Покажем, что существует обратное отображение к 
$$
\begin{cases}
\theta(t) = \arsh{\Big( \sh{(\psi(0))}\sin{t} \Big)}\\
\varphi(t) = \begin{cases}
\arctg{\Big( \ch{(\psi(0))}\tg{t} \Big)},\ t \in (0, \pi/2]\\
\pi - \arctg{\Big( \ch{(\psi(0))}[\tg{(\pi - t)}] \Big)},\ t \in (\pi/2, \pi)
\end{cases}
\end{cases}
$$

Покажем, что формулы (\ref{X}) и (\ref{Y}) задают диффеоморфизм множества $ A = \mathbb{R} \times (0,\ \pi/2)$ и $B = \lbrace (X,Y) \in \mathbb{R}^2 \ \vert \ Y > |X| \rbrace$. Найдём явно формулы для обратного отображения. Они окажутся гладкими функциями.
$$
\psi_0 \in \mathbb{R},\ t_0 \in (0,\ \pi/2)
$$
$$
\begin{cases}
X = \sh{\psi_0}\sin{t_0}\\
Y = \ch{\psi_0}\tg{t_0}
\end{cases}
$$
Выразим из первого уравнения, подставим во второе и воспользуемся формулами для тригонометрических и гиперболических функций:
$$
\sh{\psi_0} = \frac{X}{\sin{t_0}} \Rightarrow Y = \frac{\sin{t_0}}{\cos{t_0}} \sqrt{1+\frac{X^2}{\sin^2{t_0}}} = \frac{\sqrt{X^2 + \sin^2{t_0}}}{\cos{t_0}} \Rightarrow
$$
$$
\Rightarrow Y^2 = \frac{X^2 + \sin^2{t_0}}{\cos^2{t_0}} = \frac{X^2+\sin^2{t_0}}{1-\sin^2{t_0}} \Rightarrow
$$
$$
\Rightarrow Y^2(1-\sin^2{t_0}) = \sin^2{t_0}+X^2 \Leftrightarrow Y^2-X^2 = (1+Y^2)\sin^2{t_0} \Leftrightarrow \sin^2{t_0} = \frac{Y^2-X^2}{1+Y^2} 
$$
$$
\Leftrightarrow t_0 = \arcsin{\sqrt{\frac{Y^2-X^2}{1+Y^2} }} \Rightarrow \psi_0 = \arsh{\Bigg( X\sqrt{\frac{1+Y^2}{Y^2 - X^2}} \Bigg)} -\ \text{гладкие функции на множестве B} 
$$

Обратное отображение к отображению, заданному формулами (\ref{theta}), (\ref{phi}), выражается гладкими функциями:
$$
t_0 = \arcsin{\sqrt{\frac{\tg^2{\varphi}-\sh^2{\theta}}{1+\tg^2{\varphi}} }},\ \psi_0 = \arsh{\Bigg( \sh{\theta}\sqrt{\frac{1+\tg^2{\varphi}}{\tg^2{\varphi} - \sh^2{\theta}}} \Bigg)}
$$

При $t \in (0,\pi/2),\ \psi(0) \in \mathbb{R}$ мы его нашли. Теперь для $t=t_0 \in (\pi/2,\ \pi),\ \psi(0) = \psi_0 \in \mathbb{R}$, что $\lbrace (\theta, \varphi) \in \mathbb{R}^2 : \pi/2 < \varphi < \pi - \arctg{\sh{\theta}} \rbrace $.
$$
\begin{cases}
\theta(t) = \arsh{\Big( \sh{(\psi_0)}\sin{t_0} \Big)}\\
\varphi(t) = \pi - \arctg{\Big( \ch{(\psi_0)}[\tg{(\pi - t_0)}] \Big)}
\end{cases}
\Leftrightarrow
\begin{cases}
\sh{\theta} = \sh{(\psi_0)}\sin{t_0}\\
\tg{(\pi - \varphi)} = \ch{(\psi_0)}[\tg{(\pi - t_0)}]
\end{cases}
$$

Так как $t_0 \in (\pi/2,\ \pi)$, то $(\pi - t_0) \in (0,\pi/2) \Rightarrow \cos{(\pi - t_0)}=-\cos{t_0} > 0,\ \sin{(\pi - t_0)} = -\sin{t_0} > 0 \Rightarrow \tg{(\pi - t_0)} = \tg{t_0}$. По аналогичным соображениям $\tg{(\pi - \varphi)} = \tg{\varphi}$, так как $\varphi \in (\pi/2, \pi)$. Поэтому формулы выглядят так:
$$
\begin{cases}
\sh{\theta} = \sh{(\psi_0)}\sin{t_0}\\
\tg{\varphi} = \ch{(\psi_0)}\tg{t_0}
\end{cases}
$$
Это в точности формулы для случая $t_0 \in (0,\pi/2),\ \psi_0 \in \mathbb{R}$ с некоторыми изменениями, поэтому
$$
t_0 = \pi - \arcsin{\Bigg(\sqrt{\frac{\tg^2{\varphi} - \sh^2{\theta}}{1+\tg^2{\varphi}}} \Bigg)}
$$
$$
\psi_0 = \arsh{\Bigg( \sh{\theta}\sqrt{\frac{1+\tg^2{\varphi}}{\tg^2{\varphi} - \sh^2{\theta}}}\Bigg)}
$$

Надо понять, что происходит на линии $t_0 = \pi/2,\ \psi_0 \in \mathbb{R}$.
$$
\theta(\pi/2) = \arsh{\Big( \sh{(\psi_0)} \sin(\pi/2) \Big)} = \psi_0
$$
$$
\varphi(\pi/2) = \pi/2
$$

А что с формулами при $\theta = \psi_0,\ \varphi = \pi/2$:

$$
t_0 = \underset{\varphi \longrightarrow \pi/2}{\lim}\arcsin{\sqrt{\frac{\tg^2{\varphi} - \sh^2{\psi_0}}{1+\tg^2{\varphi}}}}
$$
Вычислим, воспользовавшись правилом Лопиталя:
$$
\underset{\varphi \longrightarrow \pi/2}{\lim} \frac{\tg^2{\varphi} - \sh^2{\psi_0}}{1+\tg^2{\varphi}}= (\frac{\infty}{\infty}) = \underset{\varphi \longrightarrow \pi/2}{\lim} \frac{\tg{\varphi} \frac{2}{\cos^2{\varphi}}}{\tg{\varphi}\frac{2}{\cos^2{\varphi}}} = 1
$$
Следовательно, существует предел композиции:
$$
t_0 = \arcsin{1} = \pi/2 = \pi - \arcsin{1}
$$
$$
\psi_0 = \underset{\varphi \longrightarrow \pi/2}{\lim}\arsh{\Bigg( \sh{\psi_0}\sqrt{\frac{1+\tg^2{\varphi}}{\tg^2{\varphi} - \sh^2{\psi_0}}}\Bigg)} = \arsh{\Bigg( \sh{\psi_0} \Bigg)} = \psi_0
$$
Следовательно существует диффеоморфизм областей $A'$ и $C'$.

\subsection{Множество достижимости из точки $\theta = 0,\ \varphi = 0$}

Вернёмся к системе \ref{pr1}--\ref{pr3}. Покажем, что множество достижимости это:
$$
V = \lbrace (\theta, \varphi) \in \mathbb{R}^2 : \varphi \geqslant \arctg{(\sh{\vert \theta \vert})} \rbrace
$$

Рассмотрим постоянные управления $u_1 = const, \ u_2 = const,$ такие что $u_1>0,\ -u_1 \leqslant u_2 \leqslant u_1$ и найдём соответствуюшие траектории векторного поля с начальным условием $\theta(0)=0,\ \varphi(0) = 0$:
\begin{equation}
\label{sys}
\begin{cases}
\dot{\theta} = u_2\\
\dot{\varphi} = \frac{u_1}{\ch{\theta}}
\end{cases}
\end{equation}
Решение:
$$
\begin{cases}
u_2 =0,\ \theta(t) \equiv 0,\ \varphi(t) = u_1 t\\
u_2 \neq 0,\ \theta(t) = u_2t,\ \varphi(t) = \frac{u_1}{u_2}\arctg{(\sh{(u_2t)})}
\end{cases}
$$

\begin{figure}[H]
\label{pic3}
\caption{Траектории при $u_2 = 0,$ а также $u_2 = \pm 1,\ u_1 \in \lbrace 1,\ 3,\ 5,\ 10,\ 20,\ 50  \rbrace$}
\includegraphics[width = 10 cm, height = 8 cm]{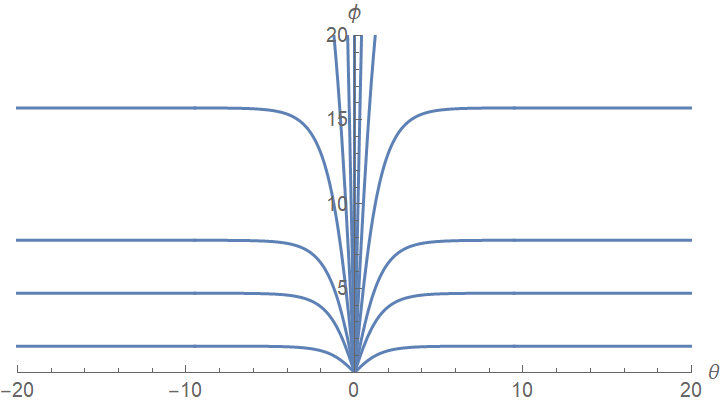}
\centering
\end{figure}

При $u_2 = 0$ получаем вертикальную прямую $\theta=0$. Она разбивает $V$ на два непересекающихся множества: $V_{+} = V \cap \lbrace \theta > 0 \rbrace,\ V_{-} = V \cap \lbrace \theta < 0 \rbrace$.\\
Пусть $u_2 = \pm 1 \Rightarrow u_1 \geqslant 1$ и $u_2 t = \pm t,\ \frac{u_1}{u_2} = \pm u_1 = u, u \geqslant 1\ \text{или} \ u \leqslant -1$. Поэтому
$$
\begin{cases}
\theta = \pm t\\
\varphi = u \arctg{(\sh{(\pm t)})} 
\end{cases}
\Rightarrow
\begin{cases}
\pm t = \theta \\
u = \frac{\varphi}{\arctg{(\sh{\theta})}}
\end{cases}
$$
Отсюда понятно, что $V_{+}$ биективно $\lbrace (s,u) : s>0,\ u \geqslant 1 \rbrace$, а $V_{-}$ биективно $\lbrace (s,u) : s<0,\ u \leqslant -1 \rbrace$. Значит, что все точки множества $V$ достижимы из $\theta = 0,\ \varphi = 0$.\\

Теперь покажем, что любое векторное поле $u_1 X_1 + u_2 X_2$ при любом допустимом $u$ в каждой точке границы $V$ касается кривой $\varphi = \arctg{(\sh{\vert \theta \vert})}$ или же направлено внутрь области $int(V)$. Для этого рассмотрим дифференциальные формы в соответствующих областях:
$$
\begin{cases}
d \varphi - \frac{d \theta}{\ch{\theta}} = d (\varphi - \arctg{(\sh{\theta})}),\ \theta > 0\\
d \varphi + \frac{d \theta}{\ch{\theta}} = d (\varphi - \arctg{(\sh{(-\theta)})}),\ \theta < 0
\end{cases}
$$

В точках границы $V \setminus (0,0)$ при подстановке в эти формы произвольного векторного поля $v$ знак результата укажет, направлено поле $v$ внутрь $V$, по касательной к границе или наружу. Подставим наше поле $u_1 X_1 + u_2 X_2 = \frac{u_1}{\ch{\theta}} \frac{\partial}{\partial \varphi} + u_2 \frac{\partial}{\partial \theta}$ при $ u_1^2 - u_2^2 \geqslant 0,\ u_1 > 0 \Rightarrow u_1 + u_2 \geqslant 0,\ u_1 - u_2 \geqslant 0$ (что следует из формулы $u_1^2-u_2^2 = (u_1-u_2)(u_1+u_2)$ и того, что $u_1 > 0$). Поэтому
$$
\begin{cases}
\theta > 0,\ \frac{u_1}{\ch{\theta}} - \frac{u_2}{\ch{\theta}} \geqslant 0\\
\theta < 0,\ \frac{u_1}{\ch{\theta}} + \frac{u_2}{\ch{\theta}} \geqslant 0
\end{cases}
$$

Отсюда следует, что множество достижимости из точки $\theta = 0,\ \varphi = 0$ за неотрицательное время равно $V$.

\subsection{Область существования оптимальных траекторий}

Сначала покажем, что множество достижимости из произвольной точки $q_0$ такое же, как и из $q_0 = 0$.\\
Решение системы (\ref{sys}) для $q_0 = (\theta_0,\varphi_0)$ и постоянных $u_1,u_2$:
$$
\begin{cases}
u_2 =0,\ \theta(t) \equiv \theta_0,\ \varphi(t) = u_1 t + \varphi_0\\
u_2 \neq 0,\ \theta(t) = u_2t + \theta_0,\ \varphi(t) = \frac{u_1}{u_2}\arctg{(\sh{(u_2t)})} + \varphi_0
\end{cases}
$$

Поэтому по аналогичным соображениям, приведённым в предыдущей секции, множество достижимости из точки $q_0$:
$$
V_{q_0} = \lbrace (\theta, \varphi) \in \mathbb{R}^2 : \varphi \geqslant \arctg{(\sh{\vert \theta - \theta_0 \vert})}+ \varphi_0 \rbrace
$$

Также стоит отметить, что при отрицательном времени 
$$
V_{q_0}^{-} = \lbrace (\theta, \varphi) \in \mathbb{R}^2 : \varphi \leqslant -\arctg{(\sh{\vert \theta - \theta_0 \vert})}+ \varphi_0 \rbrace
$$

Поэтому, доказав, что  
$$
\mathcal B = \{(\theta, \varphi) \in M \vert \pi - \arctan \sinh |\theta| >  \varphi \geqslant \arctan \sinh |\theta|\} - 
$$
множество достижимости оптимальных траекторий из $(0,0)$, то мы докажем, что для произвольного $q_0$ это будет множество
$$
\mathcal B_{q_0} = \{(\theta, \varphi) \in M \vert \varphi_0 + \pi - \arctg{ \sh{ |\theta-\theta_0|}} >  \varphi \geqslant \varphi_0 + \arctg{ \sh {|\theta-\theta_0|}}\},
$$
Рассмотрим точку $q_1 \in \mathcal B$.
Проверим, что для $q_1 = (\theta_1,\varphi_1)$ выполняются условия теоремы 2 из статьи \cite{sl_exist}:
\begin{enumerate}
\item
Первое условие выполняется, так как $\mathcal B_{q_0} \subset V_{q_0}$.
\item
Покажем, что $V_{(0,0)} \cap V_{q_1}^{-}$ - компакт. Действительно, ведь их границы пересекаются в двух точках так, что мы получаем ограничения слева и справа, а ограничения сверху и снизу заданы самими кривыми $\Phi_1 = \arctg{(\sh{|\theta}|)}$ и $\Phi_2 = \varphi_1 - \arctg{(\sh{|\theta - \theta_0|})}$. Как пересекаются эти кривые?\\
Так как обе они состоят из трёх частей, соответственно: вершина $(0,0)$ для $\Phi_1$ ($q_1$ для $\Phi2$), правая кривая при $\theta>0$ для $\Phi_1$ ($\theta - \theta_0 > 0$ для $\Phi2$), левая кривая при $\theta <0$ для $\Phi_1$ ($\theta - \theta_0 <0$ для $\Phi2$). $\underset{\theta \rightarrow \pm \infty}{\lim} \Phi_1 = \pi/2$, $\underset{\theta \rightarrow \pm \infty}{\lim} \Phi_2 = \varphi_1 - \pi/2$, причём стремление монотонное (производная имеет постоянный знак). Так как $0 \leqslant \varphi_1 < \pi$, то $-\pi/2 \leqslant \varphi_1 - \pi/2 < \pi/2$.
Это значит, что $\Phi_1 \vert_{\theta > 0} \cap \Phi_2 \vert_{\theta-\theta_0 > 0}$ в единственной точке $(\tilde{\theta}_{+},\tilde{\varphi}_{+})$ (самой крайней левой), как и $\Phi_1 \vert_{\theta < 0} \cap \Phi_2 \vert_{\theta - \theta_0 < 0}$ в единственной точке $(\tilde{\theta}_{-},\tilde{\varphi}_{-})$ (самой крайней правой). Поэтому пересечение $V_{(0,0)} \cap V_{q_1}^{-}$ можно заключить в прямоугольник $\lbrace (\theta,\varphi): \tilde{\theta}_{-} \leqslant \theta \leqslant \tilde{\theta}_{+},\ \tilde{\varphi}_{-} \leqslant \varphi \leqslant \tilde{\varphi}_{+} \rbrace$, из чего следует ограниченность.\\
Замкнутость следует из того, что границы обоих множеств включены.
\\

\item
Покажем, что 
$$
\sup \lbrace t_1 > 0: \exists \ \text{траектория}\ q(t) \ \text{системы}\ (\ref{pr1}), (\ref{pr2}),\ t \in [0,t_1]: q(0) = q_0,\ q(t_1)=q_1 \rbrace < +\infty
$$ 
Это следует из того, что $\dot{\theta} = \frac{u_1}{\ch{\theta}} > 0,$ то есть $\theta_1$ достижимо за конечное время при любых $u_1, u_2 \in U$.
\end{enumerate}

\subsection{Оптимальный синтез на нижней границе области $C'$}

Всё написано в секции 1.11.

\subsection{Оптимальный синтез в области $C'$}

Исходя из того, что для любой точки $q_1 = (\theta_1,\varphi_1) \in C'$ существует оптимальная траектория, а также из того, что экстремальная траектория, проходящая через $q_1$, единственная, следует, что для каждой точки $q_1 \in C'$ оптимальная траектория, соединяющая её с $(0,0)$, находится единственным образом из формул, написанных в предыдущей секции. Надо сначала найти $\psi_{q_1}$:
$$
\psi_{q_1} = \arsh{\Bigg( \sh{\theta_1}\sqrt{\frac{1+\tg^2{\varphi_1}}{\tg^2{\varphi_1} - \sh^2{\theta_1}}}\Bigg)}
$$
Далее, находим момент, в который мы достигаем точки $q_1$:
$$
t_{q_1} = 
\begin{cases}
\arcsin{\sqrt{\frac{\tg^2{\varphi_1} - \sh^2{\theta_1}}{1+\tg^2{\varphi_1}}}},\ \varphi_1 \in (0,\pi/2)\\
\pi/2,\ \varphi_1 = \pi/2 \\
\pi - \arcsin{\sqrt{\frac{\tg^2{\varphi_1} - \sh^2{\theta_1}}{1+\tg^2{\varphi_1}}}},\ \varphi_1 \in (\pi/2,\pi)
\end{cases}
$$

Оптимальная траектория, соединяющая точки $(0,0)$ и $q_1 = (\theta_1,\varphi_1) \in C'$.\\

Если $t_{q_1} < \pi/2$, то

\begin{equation}
\label{synthesis1}
\begin{cases}
\theta(t) = \arsh{\Big( \sh{(\psi_{q_1})}\sin{t} \Big)}\\
\varphi(t) = 
\arctg{\Big( \ch{(\psi_{q_1})}\tg{t} \Big)},\ t \in (0, t_{q_1})
\end{cases}
\end{equation}

Если $t_{q_1} = \pi/2$, то

\begin{equation}
\label{synthesis2}
\begin{cases}
\theta(t) = \arsh{\Big( \sh{(\psi_{q_1})}\sin{t} \Big)}\\
\varphi(t) = 
\begin{cases}
\arctg{\Big( \ch{(\psi_{q_1})}\tg{t} \Big)},\ t \in (0, t_{q_1})\\
\pi/2
\end{cases}
\end{cases}
\end{equation}

Если $\pi/2 < t_{q_1} < \pi$, то

\begin{equation}
\label{synthesis3}
\begin{cases}
\theta(t) = \arsh{\Big( \sh{(\psi_{q_1})}\sin{t} \Big)}\\
\varphi(t) = \begin{cases}
\arctg{\Big( \ch{(\psi_{q_1})}\tg{t} \Big)},\ t \in (0, \pi/2)\\
\pi/2,\ t = \pi/2\\
\pi - \arctg{\Big( \ch{(\psi_{q_1})}[\tg{(\pi - t)}] \Big)},\ t \in (\pi/2, t_{q_1}),
\end{cases}
\end{cases}
\end{equation}

\subsection{Программа для Wolfram Mathematica для вычисления траектории в области $C'$}

\begin{verbatim}
\[Psi][\[Theta]_, \[Phi]_] := 
  ArcSinh[Sinh[\[Theta]] Sqrt[(1 + Tan[\[Phi]]^2)/(Tan[\[Phi]]^2 - 
        Sinh[\[Theta]]^2)]];
\[Tau]1[\[Theta]_, \[Phi]_] := 
  ArcSin[Sqrt[(Tan[\[Phi]]^2 - Sinh[\[Theta]]^2)/(1 + Tan[\[Phi]]^2)]];
\[Tau]2[\[Theta]_, \[Phi]_] := \[Pi] - 
   ArcSin[Sqrt[(Tan[\[Phi]]^2 - Sinh[\[Theta]]^2)/(1 + 
        Tan[\[Phi]]^2)]];
Ph1[t_, ps_] := ArcTan[Cosh[ps] Tan[t]];
Ph2[t_, ps_] := \[Pi] - ArcTan[Cosh[ps] Tan[\[Pi] - t]];
Th[t_, ps_] := ArcSinh[Sinh[ps] Sin[t]];
Manipulate[
 If[p[[2]] <= \[Pi]/2, 
  Show[ParametricPlot[{Th[t, \[Psi][p[[1]], p[[2]]]], 
     Ph1[t, \[Psi][p[[1]], p[[2]]]]}, {t, 0, \[Tau]1[p[[1]], p[[2]]]},
     AxesLabel -> {\[Theta], \[Phi]}, 
    PlotRange -> {{-n, n}, {0, \[Pi]}}], 
   Plot[ArcTan[Sinh[t]], {t, 0, n}, PlotStyle -> Green], 
   Plot[-ArcTan[Sinh[t]], {t, -n, 0}, PlotStyle -> Green], 
   Plot[\[Pi] - ArcTan[Sinh[t]], {t, 0, n}, PlotStyle -> Green], 
   Plot[\[Pi] + ArcTan[Sinh[t]], {t, -n, 0}, PlotStyle -> Green]], 
  Show[ParametricPlot[{Th[t, \[Psi][p[[1]], p[[2]]]], 
     Ph1[t, \[Psi][p[[1]], p[[2]]]]}, {t, 0, \[Pi]/2}, 
    AxesLabel -> {\[Theta], \[Phi]}, 
    PlotRange -> {{-n, n}, {0, \[Pi]}}], 
   ParametricPlot[{Th[t, \[Psi][p[[1]], p[[2]]]], 
     Ph2[t, \[Psi][p[[1]], p[[2]]]]}, {t, \[Pi]/2, \[Tau]2[p[[1]], 
      p[[2]]]}], Plot[ArcTan[Sinh[t]], {t, 0, n}, PlotStyle -> Green],
    Plot[-ArcTan[Sinh[t]], {t, -n, 0}, PlotStyle -> Green], 
   Plot[\[Pi] - ArcTan[Sinh[t]], {t, 0, n}, PlotStyle -> Green], 
   Plot[\[Pi] + ArcTan[Sinh[t]], {t, -n, 0}, 
    PlotStyle -> Green]]], {{p, {0, \[Pi]/3}}, Locator}, {{n, 3}, 0.1,
   10000}]
\end{verbatim}

\subsection{Точки выше верхней границы $C'$}

Зададим семейство допустимых кривых, зависящих от $\alpha > 0$, соединяюших начало координат с точкой $q_1 = (\theta_1,\varphi_1)$, т.ч. $\varphi_1 > \pi - \arctg{(\sh{|\theta|})}$. Кривая состоит из 3 частей:
1) Двигаемся по правой нижней границе до точки $(\theta,\varphi) = (\alpha,\arctg{(\sh{\alpha})})$, если $\theta_1 > 0$, по левой нижней границе до точки $(-\alpha, -\arctg{(\sh{(-\alpha)})})$, если $\theta_1 < 0$;\\
2) Двигаемся вертикально вверх до пересечения с перенесённой правой верхней границей, проходящей через точку $(\theta_1,\varphi_1)$ при $\theta_1 > 0$, до пересечения с перенесённой левой верхней границе, проходящей через точку $(\theta_1,\varphi_1)$ при $\theta_1 < 0$;\\
3) Двигаемся по перенесённой верхней границе до точки $(\theta_1,\varphi_1)$.\\

Какими управлениями задаются каждая из этих частей?\\
Вспомним, как выглядит дифференциальное уравление:
$$
\begin{cases}
\dot{\theta} = u_2\\
\dot{\varphi} = \frac{u_1}{\ch{\theta}}
\end{cases}
$$
1) $u_1 = 1,\ u_2 = 1$ при $\theta_1 > 0$, $u_1 = 1,\ u_2 = -1$ при $\theta_1 <0$; $t \in [0,t_1]$;\\
2) $u_1 = 1,\ u_2 = 0$; $t \in [t_1,t_2]$;\\
3) $u_1 = 1,\ u_2 = -1$ при $\theta_1 > 0$, $u_1 = 1,\ u_2 = 1$ при $\theta_1 < 0$; $t \in [t_2,t_3]$.\\

Вычислим значение функционала длины на такой кривой:
$$
\int_0^{t_3} \sqrt{u_1^2-u_2^2} \,dt = ( \int_0^{t_1} + \int_{t_1}^{t_2} + \int_{t_2}^{t_3}) \sqrt{u_1^2-u_2^2} \,dt = \int_{t_1}^{t_2} \,dt = t_2-t_1  
$$
Поэтому нам на самом деле нужно найти $t_2$ и $t_1$.\\
Сначала для $\theta_1 > 0$:\\
1) С начальным условием $\theta(0) = 0,\ \varphi(0) = 0$ получаем решение $(\theta(t),\varphi(t)) = (t,\arctg{\sh{t}})$ на отрезке $t \in [0,\alpha]$;\\
2) С начальным условием $\theta(\alpha) = \alpha,\ \varphi(\alpha) = \arctg{(\sh{\alpha})}$ получаем решение $(\theta(t),\varphi(t)) = ( \alpha, \frac{t}{\ch{\alpha}} - \frac{\alpha}{\ch{\alpha}} + \arctg{(\sh{\alpha})} )$ на отрезке $t \in [\alpha, t_2]$;\\
3) С начальным условием $\theta(t_3) = \theta_1,\ \varphi(t_3) = \varphi_1$ получаем решение $(\theta(t),\varphi(t)) = (-t+t_3+\theta_1, -\arctg{(\sh{(-t+t_3+\theta_1)})}+\varphi_1+\arctg{(\sh{(\theta_1)})})$ на отрезке $t \in [t_2,t_3]$.\\
$t_1=\alpha$ мы нашли.\\

Осталось найти $t_2$ из пересечения вертикальной прямой 2 с кривой 3. В точке $t_2$ вертикальная прямая 2 достигает точку $-\arctg{(\sh{\alpha})}+\varphi_1+\arctg{(\sh{(\theta_1)})}$.\\
$$
\frac{t_2}{\ch{\alpha}} - \frac{\alpha}{\ch{\alpha}} + \arctg{(\sh{\alpha})} = -\arctg{(\sh{\alpha})}+\varphi_1+\arctg{(\sh{(\theta_1)})} \Leftrightarrow
$$
$$
\Leftrightarrow t_2 = [ \frac{\alpha}{\ch{\alpha}} -2\arctg{(\sh{\alpha})}+\varphi_1+\arctg{(\sh{(\theta_1)})} ]\ch{\alpha}
$$
И теперь длину кривой:
$$
t_2 - t_1 = [ \frac{\alpha}{\ch{\alpha}} -2\arctg{(\sh{\alpha})}+\varphi_1+\arctg{(\sh{(\theta_1)})} ]\ch{(\alpha)} - \alpha =
$$
$$
= [-2\arctg{(\sh{\alpha})}+\varphi_1+\arctg{(\sh{(\theta_1)})} ]\ch{(\alpha)} = L(\alpha)
$$
Вычислим предел:

$$
\underset{\alpha \rightarrow +\infty}{\lim} L(\alpha) = \underset{\alpha \rightarrow +\infty}{\lim}\Bigg( [-2\arctg{(\sh{\alpha})}+\varphi_1+\arctg{(\sh{(\theta_1)})} ]\ch{(\alpha)} \Bigg) =
$$
$$
= \infty \cdot \sign{\Big( \varphi_1+ \arctg{(\sh{\theta_1})} - \pi \Big)} = +\infty,
$$

так как $\varphi_1 > \pi - \arctg{\theta_1} \Leftrightarrow \varphi_1+ \arctg{(\sh{\theta_1})} - \pi > 0 $.

Теперь $\theta_1 < 0$:\\

1) С начальным условием $\theta(0) = 0,\ \varphi(0) = 0$ получаем решение $(\theta(t),\varphi(t)) = (-t,-\arctg{\sh{(-t)}})$ на отрезке $t \in [0,\alpha]$;\\
2) С начальным условием $\theta(\alpha) = -\alpha,\ \varphi(\alpha) = -\arctg{(\sh{(-\alpha)})}$ получаем решение $(\theta(t),\varphi(t)) = ( -\alpha, \frac{t}{\ch{(-\alpha)}} + \frac{\alpha}{\ch{(-\alpha)}} - \arctg{(\sh{(-\alpha)})} )$ на отрезке $t \in [\alpha, t_2]$;\\
3) С начальным условием $\theta(t_3) = \theta_1,\ \varphi(t_3) = \varphi_1$ получаем решение $(\theta(t),\varphi(t)) = (t-t_3+\theta_1, \arctg{(\sh{(t-t_3+\theta_1)})}+\varphi_1-\arctg{(\sh{(\theta_1)})})$ на отрезке $t \in [t_2,t_3]$.\\
$t_1=\alpha$ мы нашли.\\

Осталось найти $t_2$ из пересечения вертикальной прямой 2 с кривой 3. В точке $t_2$ вертикальная прямая 2 достигает точку $\arctg{(\sh{(-\alpha)})}+\varphi_1-\arctg{(\sh{(\theta_1)})}$.\\
$$
\frac{t_2}{\ch{(-\alpha)}} + \frac{\alpha}{\ch{(-\alpha)}} - \arctg{(\sh{(-\alpha)})} = \arctg{(\sh{(-\alpha)})}+\varphi_1-\arctg{(\sh{(\theta_1)})} \Leftrightarrow
$$
$$
\Leftrightarrow t_2 = [ -\frac{\alpha}{\ch{(-\alpha)}} + 2\arctg{(\sh{(-\alpha)})}+\varphi_1-\arctg{(\sh{(\theta_1)})} ]\ch{(-\alpha)}
$$
И теперь длину кривой:
$$
t_2 - t_1 = [ -\frac{\alpha}{\ch{(-\alpha)}} +2\arctg{(\sh{(-\alpha)})}+\varphi_1-\arctg{(\sh{(\theta_1)})} ]\ch{(-\alpha)} - \alpha =
$$
$$
= [-2\frac{\alpha}{\ch{\alpha}}-2\arctg{(\sh{\alpha})}+\varphi_1-\arctg{(\sh{(\theta_1)})} ]\ch{\alpha} = L(\alpha)
$$
Вычислим предел:

$$
\underset{\alpha \rightarrow +\infty}{\lim} L(\alpha) = \underset{\alpha \rightarrow +\infty}{\lim}\Bigg( [-2\frac{\alpha}{\ch{\alpha}}-2\arctg{(\sh{\alpha})}+\varphi_1-\arctg{(\sh{(\theta_1)})} ]\ch{(\alpha)} \Bigg) =
$$
$$
=  \infty \cdot \sign{\Big( \varphi_1- \arctg{(\sh{\theta_1})} - \pi \Big)} = +\infty,
$$

так как $\theta_1 < 0$, $\varphi_1 > \pi + \arctg{\theta_1} \Leftrightarrow \varphi_1 - \arctg{(\sh{\theta_1})} - \pi > 0 $.

Теперь разберём случай, когда $\theta_1 = 0,\ \varphi_1 > \pi$. Рассмотрим кривые, которые идут справа (как когда $\theta_1 > 0$).\\

1) С начальным условием $\theta(0) = 0,\ \varphi(0) = 0$ получаем решение $(\theta(t),\varphi(t)) = (t,\arctg{\sh{t}})$ на отрезке $t \in [0,\alpha]$;\\
2) С начальным условием $\theta(\alpha) = \alpha,\ \varphi(\alpha) = \arctg{(\sh{\alpha})}$ получаем решение $(\theta(t),\varphi(t)) = ( \alpha, \frac{t}{\ch{\alpha}} - \frac{\alpha}{\ch{\alpha}} + \arctg{(\sh{\alpha})} )$ на отрезке $t \in [\alpha, t_2]$;\\
3) С начальным условием $\theta(t_3) = 0,\ \varphi(t_3) = \varphi_1$ получаем решение $(\theta(t),\varphi(t)) = (-t+t_3, -\arctg{(\sh{(-t+t_3)})}+\varphi_1)$ на отрезке $t \in [t_2,t_3]$.\\
$t_1=\alpha$ мы нашли.\\

Осталось найти $t_2$ из пересечения вертикальной прямой 2 с кривой 3. В точке $t_2$ вертикальная прямая 2 достигает точку $-\arctg{(\sh{\alpha})}+\varphi_1+\arctg{(\sh{(\theta_1)})}$.\\
$$
\frac{t_2}{\ch{\alpha}} - \frac{\alpha}{\ch{\alpha}} + \arctg{(\sh{\alpha})} = -\arctg{(\sh{\alpha})}+\varphi_1 \Leftrightarrow
$$
$$
\Leftrightarrow t_2 = [ \frac{\alpha}{\ch{\alpha}} -2\arctg{(\sh{\alpha})}+\varphi_1 ]\ch{\alpha}
$$
И теперь длину кривой:
$$
t_2 - t_1 = [ \frac{\alpha}{\ch{\alpha}} -2\arctg{(\sh{\alpha})}+\varphi_1 ]\ch{(\alpha)} - \alpha =
$$
$$
= [-2\arctg{(\sh{\alpha})}+\varphi_1+\arctg{(\sh{(\theta_1)})} ]\ch{(\alpha)} = L(\alpha)
$$
Вычислим предел:

$$
\underset{\alpha \rightarrow +\infty}{\lim} L(\alpha) = \underset{\alpha \rightarrow +\infty}{\lim}\Bigg( [-2\arctg{(\sh{\alpha})}+\varphi_1]\ch{(\alpha)} \Bigg) =
$$
$$
=  \infty \cdot \sign{\Big( \varphi_1- \pi \Big)} = +\infty,
$$

так как $\varphi_1 > \pi$.

\subsection{Точки верхней границы $C'$}

Предположим, что существует оптимальная траектория, соединяющая $q_0$ и $q_1 = (\theta_1, \pi- \arctg{(\sh{|\theta_1|})})$, $\theta_1 \neq 0$. Тогда она удовлетворяет Принципу максимума Понтрягина. Однако мы нашли все траектории, удовлетворяющие ПМП. Анормальные траектории задают нижнюю границу множества $C'$, а все нормальные траектории задаются формулами (\ref{full_resh}). При $n=1$ эти формулы выглядят так:
$$
\begin{cases}
\theta(t) = \arsh{\Big( \sh{(\psi_0)}\sin{t} \Big)}\\
\varphi(t) = \begin{cases}
\arctg{\Big( \ch{(\psi_0)}\tg{t} \Big)},\ t \in (0, \pi/2]\\
\pi - \arctg{\Big( \ch{(\psi_0)}[\tg{(\pi - t)}] \Big)},\ t \in (\pi/2, \pi]
\end{cases}
\end{cases}
$$
Мы показали, что $(0,\pi)_t \times \mathbb{R}_{\psi_0}$ диффеоморфно области $C'$, а при $t = 0$ мы находимся в точке $q_0$, в то время как при $t=\pi$ мы попадаем в точку $(0,\pi)$ для любого $\psi_0 \in \mathbb{R}$. Поэтому для точек верхней не существует оптимальных траекторий. Для них даже не существует экстремальных траекторий (удовлетворяющих ПМП).\\
Если $\tilde{q} = (0,\pi)$, то для неё существует континуум оптимальных траекторий, соединяющих $q_0$ и $\tilde{q}$, что видно из формул.\\

Лемма 4.4 из \cite{beem}:\
{\it Для Лоренцева расстояния $d$, если $d(p,q) < \infty$, $p_n \rightarrow p$, и $q_n \rightarrow q$, то $d(p,q) \leqslant \liminf{d(p_n,q_n)}$. Если $d(p,q) = \infty$, $p_n \rightarrow p$, и $q_n \rightarrow q$, то $\underset{n \rightarrow \infty}{\lim} d(p_n,q_n) = \infty$.}\\
Так как наша задача является Лоренцевой, то мы можем воспользоваться этой леммой.\\
Рассмотрим последовательность точек $q_n \rightarrow \tilde{q}$, где $q_n \in C'$, а $\tilde{q} \in \lbrace (\theta,\varphi) : \theta \in \mathbb{R},\ \varphi = \pi - \arctg{( \sh{|\theta|} )} \rbrace$. Предположим, что $\underset{n \rightarrow \infty}{\lim} d(p,q) = \infty$. Тогда, по лемме, $\underset{n \rightarrow \infty}{\lim} d(q_0,q_n) = \infty$. Но длину кривой выражает время $t$, а по нашим вычислениям, для всех $q_n$ существует оптимальная траектория, и мы показали, что $0 < t < \pi $. Значит, что $\underset{n \rightarrow \infty}{\lim} d(q_0,q_n) \leqslant \pi$. Поэтому $d(q_0,\tilde{q}) \leqslant \pi$.\\
Теперь приведём последовательность $q_n \in C',\ q_n \rightarrow \tilde{q} = (\tilde{\theta}, \pi - \arctg{(\sh{|\tilde{\theta}|})})$: $q_n = (\theta_n,\varphi_n) = (\tilde{\theta},\pi - \arctg{(\sh{|\tilde{\theta}|})} - \frac{1}{n}),\ n \in \mathbb{N},\ n \geqslant n_0$, выбрав начальный $n_0$ таким, что $q_n \in C'$, например, из условия, что $\pi - \arctg{(\sh{|\tilde{\theta}|})} - \frac{1}{n_0} \geqslant \pi/2 \Leftrightarrow n_0 ( \pi/2 -\arctg{(\sh{|\tilde{\theta}|})} ) \geqslant 1 \Leftrightarrow n_0 \geqslant \frac{1}{\pi/2 -\arctg{(\sh{|\tilde{\theta}|})} } $. Вычислим
$$
d(q_0,q_n) = \pi - \arcsin{\sqrt{\frac{\tg^2{\varphi_n} - \sh^2{\theta_n}}{1+\tg^2{\varphi_n}}}} = \pi - \arcsin{\sqrt{\frac{\tg^2{[\arctg{(\sh{|\tilde{\theta}|})} + \frac{1}{n}]} - \sh^2{\tilde{\theta}}}{1+\tg^2{[\arctg{(\sh{|\tilde{\theta}|})} + \frac{1}{n}]}}}}
$$
$$
\underset{n \rightarrow \infty}{\lim} d(q_0,q_n) = \pi
$$
Поэтому $\forall \tilde{q} \in \lbrace (\theta,\varphi) : \theta \in \mathbb{R},\ \varphi = \pi - \arctg{( \sh{|\theta|} )} \rbrace$ имеем $ d(q_0,\tilde{q}) = \pi$.

\subsection{Поля Киллинга}

По определению, векторное поле называется полем Киллинга если, производная Ли метрики вдоль него равно нулю. Постараемся найти эти поля, воспользовавшись равенством (11.1) из \cite{lor_lob}. Если $X$ - поле Киллинга, $V,\ W$ - векторные поля, то

$$
X(g(V,W)) = g([X,V],W)+g(V,[X,W])
$$

Распишем искомые поля по базису из собственных векторов метрики $g =  d \theta^2 - \ch^2{\theta}d \varphi^2$, $X = c_1 X_1 + c_2 X_2 = c_1 \frac{1}{\ch{\theta}} \frac{\partial}{\partial \varphi}+ c_2\frac{\partial }{\partial \theta}$

Составим уравнения на коэффициенты (сначала для $V=X_1,\ W = X_2$):

$$
[X,X_1] = [c_1 X_1+c_2X_2,X_1] = [\frac{c_1}{\ch{\theta}} \frac{\partial}{\partial \varphi}+c_2 \frac{\partial}{\partial \theta},\frac{1}{\ch{\theta}}\frac{\partial}{\partial \varphi}] =
$$
$$
= ( c_2 \partial_{\theta}(\frac{1}{\ch{\theta}})-\frac{1}{\ch{\theta}}\partial_{\varphi}(\frac{c_1}{\ch{\theta}}) )\frac{\partial}{\partial \varphi} + ( -\frac{1}{\ch{\theta}}\partial_{\varphi}(c_2) )\frac{\partial}{\partial \theta}=
$$
$$
= -\frac{c_2 \sh{\theta}}{\ch^2{\theta}}\frac{\partial}{\partial \varphi}-\frac{\partial_{\varphi}c_1}{\ch^2{\theta}}\frac{\partial}{\partial \varphi}-\frac{\partial_{\varphi}c_2}{\ch{\theta}}\frac{\partial}{\partial \theta} =-\frac{1}{\ch{\theta}}(c_2\sh{\theta}+\partial_{\varphi}c_1)X_1 -\frac{\partial_{\varphi}c_2}{\ch{\theta}}X_2
$$
$$
[X,X_2] = [\frac{c_1}{\ch{\theta}} \frac{\partial}{\partial \varphi}+c_2 \frac{\partial}{\partial \theta},\frac{\partial}{\partial \theta}]=
$$
$$
=(-\frac{\partial_{\theta}c_1}{\ch{\theta}}+c_1\frac{\th{\theta}}{\ch{\theta}})\frac{\partial}{\partial \varphi} - \partial_{\theta}c_2 \frac{\partial}{\partial \theta} = (-\partial_{\theta}c_1+c_1\th{\theta})X_1-\partial_{\theta}c_2X_2
$$

$$
X(g(X_1,X_1)) = g([X,X_1],X_1)+g(X_1,[X,X_1]) \Leftrightarrow 0 = 2g([X,X_1],X_1) = 
$$
$$
=2g(-\frac{1}{\ch{\theta}}(c_2\sh{\theta}+\partial_{\varphi}c_1)X_1 -\frac{\partial_{\varphi}c_2}{\ch{\theta}}X_2,X_1) =\frac{2}{\ch{\theta}}(c_2\sh{\theta}+\partial_{\varphi}c_1) 
$$

$$
X(g(X_2,X_2)) = g([X,X_2],X_2)+g(X_2,[X,X_2]) \Leftrightarrow 0 = 2g([X,X_2],X_2) =
$$
$$
= 2((-\partial_{\theta}c_1+c_1\th{\theta})X_1-\partial_{\theta}c_2X_2, X_2)=-2\partial_{\theta}c_2
$$
$$
X(g(X_1,X_2)) = g([X,X_1],X_2)+g(X_1,[X,X_2]) \Leftrightarrow
$$
$$
\Leftrightarrow 0 = g(-\frac{1}{\ch{\theta}}(c_2\sh{\theta}+\partial_{\varphi}c_1)X_1 -\frac{\partial_{\varphi}c_2}{\ch{\theta}}X_2,X_2)+g(X_1,(-\partial_{\theta}c_1+c_1\th{\theta})X_1-\partial_{\theta}c_2X_2) \Leftrightarrow
$$
$$
\Leftrightarrow 0 = -\frac{\partial_{\varphi}c_2}{\ch{\theta}}-(-\partial_{\theta}c_1+c_1\th{\theta})
$$

$$
\begin{cases}
c_2\sh{\theta}+\partial_{\varphi}c_1  = 0\\
\partial_{\theta}c_2 = 0\\
\partial_{\varphi}c_2 = \ch{\theta}\partial_{\theta}c_1 - c_1\sh{\theta}
\end{cases}
$$
Из второго уравнения следует, что $c_2(\varphi,\theta) = c_2(\varphi)$. Далее, можем проинтегрировать первое уравнение
$$
\partial_{\varphi}c_1 = -c_2(\varphi)\sh{\theta} \Leftrightarrow c_1 = -\sh{\theta}\int_{0}^{\varphi}c_2(s)\,ds + f(\theta) = -\sh{\theta}u(\varphi)+f(\theta)
$$
Подставим в третье уравнение:
$$
u''(\varphi)=\ch{\theta}[-\ch{\theta}u(\varphi)+f'(\theta)]-\sh{\theta}[-\sh{\theta}u(\varphi)+f(\theta)] \Leftrightarrow
$$
$$
\Leftrightarrow u''(\varphi) = [-\ch^2{\theta}+\sh^2{\theta}]u(\varphi)+\ch{\theta}f'(\theta)-\sh{\theta}f(\theta) \Leftrightarrow
$$
$$
\Leftrightarrow u''(\varphi)+u(\varphi) = \ch{\theta}f'(\theta)-\sh{\theta}f(\theta)
$$
Видим, что левая часть зависит только от $\varphi$, а правая часть зависит только от $\theta$. Значит, что левая и правая части равны константе. Мы получаем два уравнения: на $u(\varphi)$ и на $f(\theta)$:
$$
u''(\varphi)+u(\varphi) = A = \ch{\theta}f'(\theta)-\sh{\theta}f(\theta)
$$
Оба эти уравнения линейные. Решение первого находится практически моментально:
$$
u(\varphi) = B_1\cos{\varphi}+B_2\sin{\varphi}+A
$$
Второе решается классическим методом вариации постоянной:
$$
\ch{\theta}f'(\theta) = f(\theta)\sh{\theta}
$$
$$
\frac{df}{f} = \th{\theta}d\theta
$$
$$
f = \hat{B}\ch{\theta} 
$$
Теперь варьируем постоянную $\hat{B}$:
$$
f'=\hat{B}'\ch{\theta}+\hat{B}\sh{\theta}
$$
$$
\ch{\theta}(\hat{B}'\ch{\theta}+\hat{B}\sh{\theta})-\sh{\theta}\hat{B}\ch{\theta} =A
$$
$$
\hat{B}' \ch^2{\theta} = A \Leftrightarrow \hat{B} = A\th{\theta}+B
$$
$$
f(\theta) = (A\th{\theta}+B)\ch{\theta} = A\sh{\theta}+B\ch{\theta}
$$
Таким образом коэффициенты выглядят вот так:
$$
c_1(\varphi,\theta) = -\sh{\theta}(B_1\cos{\varphi}+B_2\sin{\varphi}+A)+A\sh{\theta}+B\ch{\theta} = -\sh{\theta}(B_1\cos{\varphi}+B_2\sin{\varphi})+B\ch{\theta}
$$
$$
c_2 = u'(\varphi) = -B_1\sin{\varphi}+B_2\cos{\varphi}
$$
И мы получаем поля Киллинга:
$$
c_1X_1+c_2X_2 = (-\sh{\theta}(B_1\cos{\varphi}+B_2\sin{\varphi})+B\ch{\theta})X_1+(-B_1\sin{\varphi}+B_2\cos{\varphi})X_2
$$
3 вектора:
$$
\hat{X}_1 = \ch{\theta}X_1 = \partial_{\varphi},\ \hat{X}_2 = \sh{\theta}\cos{\varphi}X_1+\sin{\varphi}X_2,\ \hat{X}_3 = -\sh{\theta}\sin{\varphi}X_1+\cos{\varphi}X_2 
$$

Вычислим коммутаторы:
$$
[\hat{X}_1,\hat{X}_2] = [\partial_{\varphi},\frac{\sh{\theta}\cos{\varphi}}{\ch{\theta}}\partial_{\varphi}+\sin{\varphi}\partial_{\theta}] = -\frac{\sh{\theta}\sin{\varphi}}{\ch{\theta}}\partial_{\varphi}+\cos{\varphi}\partial_{\theta} =\hat{X}_3
$$
$$
[\hat{X}_1,\hat{X}_3]= [\partial_{\varphi},-\frac{\sh{\theta}\sin{\varphi}}{\ch{\theta}}\partial_{\varphi}+\cos{\varphi}\partial_{\theta}]=\frac{\sh{\theta}\cos{\varphi}}{\ch{\theta}}\partial_{\varphi}-\sin{\varphi}\partial_{\theta} =- \hat{X}_2
$$
$$
[\hat{X}_2,\hat{X}_3]=[\frac{\sh{\theta}\cos{\varphi}}{\ch{\theta}}\partial_{\varphi}+\sin{\varphi}\partial_{\theta},-\frac{\sh{\theta}\sin{\varphi}}{\ch{\theta}}\partial_{\varphi}+\cos{\varphi}\partial_{\theta}]=
$$
$$
= \Bigg( \th{\theta}\cos{\varphi}\th{\theta}(-\cos{\varphi})+\sin{\varphi}(-\sin{\varphi})\frac{1}{\ch^2{\theta}}+\th{\theta}\sin{\varphi}\th{\theta}(-\sin{\varphi})-\cos{\varphi}\cos{\varphi} \frac{1}{\ch^2{\theta}}\Bigg)\partial_{\varphi}+
$$
$$
+\Bigg( \th{\theta}\cos{\varphi}(-\sin{\varphi})-\th{\theta}(-\sin{\varphi})\cos{\varphi} \Bigg)\partial_{\theta}=
$$
$$
= \Bigg(- \th^2{\theta}(\cos^2{\varphi}+\sin^2{\varphi}) -\frac{1}{\ch^2{\theta}}(\sin^2{\varphi}+\cos^2{\varphi}) \Bigg) \partial_{\varphi} = -\frac{\sh^2{\theta}+1}{\ch^2{\theta}}\partial_{\varphi} = -\partial_{\varphi} = -\hat{X}_1
$$
Получили такие соотношения: 
$$
[\hat{X}_1,\hat{X}_2] = \hat{X}_3,\ [\hat{X}_2,\hat{X}_3]=-\hat{X}_1,\ [\hat{X}_3,\hat{X}_1]=\hat{X}_2
$$
Поэтому алгебра Киллинга изоморфна $\mathfrak{sl}(2)$.

\subsection{Фазовые портреты полей Киллинга}

Портреты получены с помощью Wolfram Mathematica с использованием функции StreamPlot.

\subsection{Как посчитать лоренцево расстояние между произвольными двумя точками $\widetilde{H_1^2}$?}

Отметим, что поле $\hat{X}_1 = \partial_{\varphi}$ позволяет двигаться вверх-вниз.\\

Заметим также, что при $\varphi = 0$ у поля $\hat{X}_3$ первая координата зануляется, а вторая равна $1$, что позволяет двигаться вправо-влево на этой прямой. Поэтому наш маршрут будет выглядеть так:
\begin{enumerate}
\item
$\theta_0 < 0,\ \varphi_0 < 0$

Двигаемся вдоль поля $\hat{X}_1$ до пересечения с прямой $\varphi = 0$. Двигаемся вдоль поля $\hat{X}_3$ до точки $(0,0)$.

\item
$\theta_0 < 0,\ \varphi_0 \geqslant 0$

Двигаемся вдоль поля $-\hat{X}_1$ до пересечения с прямой $\varphi = 0$. Двигаемся вдоль поля $\hat{X}_3$ до точки $(0,0)$.

\item
$\theta_0 \geqslant 0,\ \varphi_0 < 0$

Двигаемся вдоль поля $\hat{X}_1$ до пересечения с прямой $\varphi = 0$. Двигаемся вдоль поля $-\hat{X}_3$ до точки $(0,0)$.

\item
$\theta_0 \geqslant 0,\ \varphi_0 \geqslant 0$

Двигаемся вдоль поля $-\hat{X}_1$ до пересечения с прямой $\varphi = 0$. Двигаемся вдоль поля $-\hat{X}_3$ до точки $(0,0)$.

\end{enumerate}


\begin{thebibliography}{99}

\bibitem{beem}
Beem, J.K., Ehrlich, P.E., Easley, K.L.: {\em Global Lorentzian Geometry}. Monographs
Textbooks Pure Appl. Math. 202, Marcel Dekker Inc. (1996)
\bibitem{notes}
А.А.~Аграчев, Ю. Л. Сачков,  
{\em Геометрическая теория управления},
Физматлит, 2005.
\bibitem{intro}
Сачков Ю.Л. {\em Введение в геометрическую теорию управления}, М.: URSS, 2021.
\bibitem{bibik}
\label{Bibik} 
Бибиков Ю. Н., {\em Курс обыкновенных дифференциальных уравнений: Учеб. пособие для ун-тов.} --- М.: Высш. шк., 1991.--- 303 с.: ил.
\bibitem{sl_exist}
Ю. Л. Сачков, Существование сублоренцевых длиннейших, {\em Дифференциальные уравнения}. 59, 12, 1702--1709 (2023)
\bibitem{lor_lob} Yu.L. Sachkov, Lorentzian distance on the Lobachevsky plane, {\em submitted}.
\bibitem{jakobson}  N. Jacobson, Lie Algebras, Dover Publications, New York, 1979
\end{thebibliography}
\end{document}